\documentclass[article,onefignum,onetabnum]{siamart171218}
\usepackage{amsfonts}
\usepackage{graphicx}
\usepackage{placeins}
\usepackage{comment}
\usepackage{MnSymbol}
\usepackage{enumerate}
\usepackage{scalerel}[2016-12-29]

\def\scaleint#1{\vcenter{\hbox{\scaleto[3ex]{\displaystyle\int}{#1}}}}
\def\bs{\mkern-15mu}
\usepackage{epstopdf}
\usepackage{epsfig}
\usepackage{wrapfig}
\usepackage{float}
\usepackage{subfigure}

\newcommand{\Int}{\scaleint{7ex}}

\newcommand{\bn}{\mathbf{n}}

\newcommand{\e}{\varepsilon}
\newcommand{\ov}{\overline}
\newcommand{\R}{\mathbb R}

\newcommand{\rr}{{\bf r}}

\newcommand{\N}{{\bf N}}

\newcommand{\beq}{\begin{equation}}
\newcommand{\eeq}{\end{equation}}

\newsiamremark{remark}{Remark}
\newsiamremark{hypothesis}{Hypothesis}
\crefname{hypothesis}{Hypothesis}{Hypotheses}
\newsiamthm{claim}{Claim}

\headers{Phase field model for climb and self-climb of dislocation loops}{Niu, Xiang, Yan}

\title{A phase field model for the motion of prismatic dislocation loops by both climb and self-climb}

\author{Xiaohua Niu\thanks{School of Mathematics and Statistics, Xiamen University of Technology, Xiamen, China
  (\email{xhniu@xmut.edu.cn; xniu@connect.ust.hk}).}
\and Xiaodong Yan\thanks{Department of Mathematics, The University of Connecticut, Storrs, CT, USA
  (\email{xiaodong.yan@uconn.edu}).}
}

\begin{document}

\maketitle

\begin{abstract}
We study the  sharp interface limit and  well-posedness of a phase field model for self-climb of prismatic dislocation loops in periodic settings. The model is set up in a  Cahn-Hilliard/Allen-Cahn framework featured with degenerate phase-dependent diffusion mobility with an additional stablizing function. Moreover, a nonlocal climb force is added to the chemical potential. We introduce a notion of weak solutions for the nonlinear model. The existence result is obtained by approximations  of the proposed model with nondegenerate mobilities. Lastly, the numerical simulations are performed to validate the phase field model and the simulation results show the big difference for the prismatic dislocation loops in the evolution time and the pattern with and without self-climb contribution.
\end{abstract}


\section{Introduction}\label{sec:intro}

In this paper, we present a phase field model for the motion of prismatic dislocation loops by both climb and self-climb.

The self-climb of dislocations (line defects in crystalline materials)  plays important roles in the properties of irradiated materials~\cite{Hirth-Lothe}. The self-climb motion
is driven by pipe diffusion of vacancies along the dislocations, and is the dominant mechanism of prismatic loop
motion and coalescence at not very high temperatures \cite{Kroupa1961,Johnson1960,DudarevReview,Swinburne2016,Okita2016,Niu2017,Niu2019,Gu2018,LiuJMPS2020,Niu2021,Dudarev2022}.
Dislocations climb is the motion of dislocations out of their slip planes with the assistance of vacancy diffusion over the bulk of the materials, and  it is an important mechanism in the plastic properties of crystalline materials at high temperatures (e.g., in dislocation creep) \cite{Hirth-Lothe,Ghoniem,Xiang2003,Xiang2006,Arsenlis,Mordehai,Keralavarma,Gu2015,Geslin2015,Zhu2018,Dudarev2022}. Phase field models (e.g., of the Cahn-Hilliard type \cite{CahHil58} or the Allen-Cahn type \cite{Allen-Cahn}) have the advantages of being able to handle topological
changes of the interfaces automatically with simple numerical implementation  on a uniform mesh of the simulation domain.  We have proposed a phase field model for self-climb of prismatic dislocation loops~\cite{Niu2021}. There are also phase field models for dislocation climb coupled with vacancy bulk diffusion~\cite{Shenyang2011,Ke-Yunzhi2014,Geslin2014,Geslin2015}. Recently, the importance of the cooperative effects of the two mechanisms of self-climb by vacancy pipe diffusion and climb by vacancy bulk diffusion has been realized~\cite{Gu2018,Dudarev2022,Kohnert2022}. To the best of our knowledge, a phase field model that accounts for the combined effect of these two mechanisms is still not available in the literature.


We propose the following modified Cahn-Hilliard Allen-Cahn type equation to model the motion of prismatic dislocation loops by both climb and self-climb:
\begin{eqnarray}\label{eqn:model0}
\label {eqn:ch} g(u)(\partial_t u+\beta \mu)&=&\nabla \cdot (M(u)\nabla \frac{\mu}{g(u)}) \text{ for } x \in \Omega \subset \mathbb{R}^2, t\in [0,\infty), \\
\label{eqn:chem} \mu&=&-\Delta u+\frac{1}{\e^2}q'(u)+\frac{1}{\e}h(u)f_{cl}.
\end{eqnarray}
In this model, without the $\beta \mu$ term on the left-hand side, it describes the self-climb of prismatic dislocation loops, and the dislocation climb by vacancy bulk diffusion is incorporated by the $\beta \mu$ term.
Here $\beta>0$ is a constant that enables a correct dislocation climb velocity, $M(u)=M_0(1-u^2)^2$, $M_0>0$, is the diffusion mobility, $q(u)=2(1-u^2)^2$ is the double well potential which takes minimums at $\pm 1$,   $g(u)=(1-u^2)^2$ is the stabilizing function which guarantees correct asymptotics in the sharp interface limit, and $\e$ is a small parameter controlling the width of the dislocation core.

In this model, we assume prismatic dislocation loops lie and evolve by self-climb in the $xy$ plane and all dislocation loops have the same perpendicular Burgers vector $\mathbb{b}=(0,0,b)$. The local dislocation line direction is given by $(\mathbf{b}/b)\times(\nabla u/|\nabla u|)$.
The last contribution $f_{cl}$ in the chemical potential $\mu$ in \eqref{eqn:chem} is the total climb force, with
\begin{equation}
f_{cl}=f_{cl}^d+f_{cl}^{app},
\end{equation}
where  $f_{cl}^d$ is the climb force generated by all the dislocations:
\begin{equation}
f_{cl}^d(x,y,u)=\frac{G b^2}{8\pi(1-\nu)}\Int_{\bs \Omega}\left(\frac{x-\ov x}{R^3}u_{\ov x}+\frac{y-\ov{y}}{R^3}u_{\ov{y}}\right)d\ov x d\ov y \label{equ:clf}
\end{equation}
 with $G$ being the shear modulus, $\nu$ the Poisson ratio, and $R=\sqrt{(x-\ov x)^2+(y-\ov y)^2}$,
and
$f_{cl}^{app}$ is the applied force.
The smooth cutoff factor $h(u)=H_0(1-u^2)^2$  is  to guarantee the climb force acts only on the dislocations. The constant $H_0>0$  is chosen such that the phase field model generates accurate climb force of the dislocations \cite{Niu2021} (c.f. Sec.~\ref{sec:asym}).

The chemical potential $\mu$  comes from variations of the classical Cahn-Hilliard energy and the elastic energy due to dislocations, i.e.
\begin{equation}
\mu=\frac{\delta E_{CH}}{\delta u}+\frac{1}{\varepsilon}h(u)\frac{\delta E_{el}}{\delta u},
\end{equation}
where
\begin{eqnarray}
&&E_{CH}(u)=\scaleint{7ex}_{\bs\Omega}\left (\frac{1}{2}|\nabla u|^2+q(u)\right) dx, \\
&& E_{el}=\scaleint{7ex}_{\bs\Omega}\left ( \frac{1}{2}u f^{d}_{cl}+u f_{cl}^{app}\right ) dx
\end{eqnarray}
are the classical Cahn-Hilliard energy and elastic energy respectively.
The climb force generated by the dislocations can be expressed as
\begin{eqnarray}
f^{d}_{cl}(x,y,u)=\frac{G b^2}{4(1-\nu)}(-\Delta)^{\frac{1}{2}}u.
\end{eqnarray}
Here $(-\Delta)^{s}u$ is  the   fractional operator defined by $$\mathcal{F}((-\Delta)^s f)=(\xi_1^2+\xi_2^2)^{\frac{s}{2}}\mathcal{F}(f)(\xi),$$ where $\xi$ is the frequency. 

This model is obtained by incorporating the dislocation climb motion into our phase field model for the self-climb motion of prismatic dislocation loops that we proposed earlier \cite{Niu2021} without the factor $g(u)$ on the left-hand side in Eq.~\eqref{eqn:model0}. This factor is mainly for the wellposedness proof, and without it, the results of dislocation velocity given by the sharp interface limit (see the remark at the end of Sec.~\ref{sec:asym}) and numerical simulations are similar.
When $g\equiv 1$ and the climb force $f_{cl}$ 
is omitted, the model reduces to the  Cahn-Hilliard/Allen-Cahn equation with degenerate mobility. Such models have attracted lots of attentions in recent years \cite{ZhaCha16}.

In this paper, we are interested in  the sharp interface limit and well-posedness  for \eqref{eqn:ch}-\eqref{eqn:chem}. Numerical simulations are also performed using the obtained phase field model.

We first derive  a sharp interface limit equation for \eqref{eqn:ch} and \eqref{eqn:chem} via formal asymptotic analysis. The following sharp interface equation is obtained as  $\e \rightarrow \infty$,
\begin{equation}
v=-\lambda\partial_{ss}\left(\alpha \kappa-H_0 f^{(0)}_{cl}(s)\right)+\eta\left(\alpha \kappa-H_0 f^{(0)}_{cl}(s)\right).
\end{equation}
 Here $\lambda$, $\alpha$ and  $\eta$ are positive  constants whose exact forms can be found in section \ref{sec:asym}.

For well-posedness of  \eqref{eqn:ch}-\eqref{eqn:chem},
 we consider the following modified problem  in a periodic setting in general dimensions. Set $\Omega=[0,2\pi]^n$, we consider
\begin{eqnarray}
\label{eqn:u} g(u)(\partial_t u+\beta \mu)&=&\nabla \cdot \left(M(u)\nabla \frac{\mu}{g(u)}\right), \hspace{0.4in} \text{ for } x \in \Omega, t\in [0,\infty) \\
\label{eqn:mu} \mu&=&-\Delta u+q'(u)+(-\Delta)^{\frac{1}{2}}u.
\end{eqnarray}
Here $g(u)=|1-u^2|^m$ for $2\leq m<\infty $, $M(u)=M_0g(u)$ for some constant $M_0>0$, $q(u)\in C^2(\R,\R)$ and there exist constants $C_i>0$, $i=1,\cdots,10$, and  $1\leq r <\infty$, such that for all $u \in \R$, \begin{eqnarray}\label{eqn:qu}C_1|u|^{r+1}-C_2\leq q(u) &\leq &C_3|u|^{r+1}+C_4, \\ \label{eqn:q'u}|q'(u)|&\leq &C_5|u|^r+C_6\\ \label{eqn:q''u}C_t|u|^{r-1}-C_8&\leq& q''(u) \leq C_9|u|^{r-1}+C_{10}.\end{eqnarray}
We see that the classical double well potential $q(u)=(1-u^2)^2$ satisfies \eqref{eqn:qu}-\eqref{eqn:q''u} with $r=3$.

In the proof, we consider  approximations  of the proposed model \eqref{eqn:u}-\eqref{eqn:mu} with positive mobilities. Given any $\theta >0$, we define
\begin{equation}
\label{eqn:gtheta} g_{\theta}(u):= \left\{ \begin{array}{cl}
                                 |1-u^2|^m & \text{  if  }|1-u^2|> \theta,\\
                                  \theta^m &\text { if } |1-u^2|\leq \theta,
                                 \end{array} \right.
\end{equation}
 and
\begin{equation}
\label{eqn:Mtheta} M_{\theta}(u):= M_0 g_{\theta}(u).
\end{equation}

Our first step is to find a sufficiently regular solution for \eqref{eqn:u}-\eqref{eqn:mu} with mobility $M_{\theta}(u)$ and stablizing function $g_{\theta}(u)$ together with a smooth potential $q(u)$. Here and throughout the paper, we use notation $\Omega_T=[0,T]\times \Omega$. This result is summarized in the following Proposition.
\begin{proposition}\label{thm:ndg}
Let  $M_{\theta},  g_{\theta}$ be defined by \eqref{eqn:Mtheta} and \eqref{eqn:gtheta}, under the assumptions \eqref{eqn:qu}-\eqref{eqn:q''u}, for any  $u_0 \in H^1(\Omega)$ and  any $T>0$, there exists a function $u_{\theta}$ such that
\begin{itemize}
\item[a)] $u_{\theta}\in L^{\infty}(0,T;H^1(\Omega))\cap C([0,T];L^p(\Omega))\cap L^2(0,T;W^{3,s}(\Omega))$, where $1\leq p <\infty$, $1\leq s <2$,
\item[b)]$\partial_tu_{\theta}\in L^2(0,T;(W^{1,q}(\Omega))')$ for $q>2$,
\item[c)]$u_{\theta}(x,0)=u_0(x) $ for all $x \in \Omega$,
\end{itemize}
which satisfies \eqref{eqn:u}-\eqref{eqn:mu} in the following weak sense
\begin{eqnarray}\notag
&&\int_0^T<\partial_tu_{\theta},\phi>_{((W^{1,q}(\Omega))',W^{1,q}(\Omega))} dt\\ \notag
&=&-\Int_{\bs \Omega_T}M_{\theta}(u_{\theta})\nabla \frac{-\Delta u_{\theta}+q'(u_{\theta})+(-\Delta)^{\frac{1}{2}}u_{\theta}}{g_{\theta}(u_{\theta})}\cdot \nabla \frac{\phi}{g_{\theta}(u_{\theta})}dxdt\\
&&-\Int_{\bs \Omega_T}\beta (-\Delta u_{\theta}+q'(u_{\theta})+(-\Delta)^{\frac{1}{2}}u_{\theta})\phi dxdt\
\label{eqn:thetandg}
\end{eqnarray}
for all $\phi\in L^2(0,T;W^{1,q}(\Omega))$ with $q>2$. In addition, the following energy inequality holds for all $t>0$.
\begin{eqnarray}
\notag
&&\scaleint{7ex}_{\bs \Omega}\left (\frac{1}{2}|\nabla u_{\theta}(x,t)|^2+q(u_{\theta}(x,t))+u_\theta(x,t)(-\Delta)^{\frac{1}{2}}u_\theta\right) dx\\ \label{eqn:engineq}
&&+\scaleint{7ex}_{\bs 0}^t\scaleint{7ex}_{\bs \Omega}M_{\theta}(u_{\theta}(x,\tau)\left|\nabla \frac{-\Delta u_{\theta}(x,\tau)+q'(u_{\theta}(x,\tau))+(-\Delta)^{\frac{1}{2}}u_{\theta}(x,\tau)}{g_{\theta}(u_{\theta}(x,\tau))}\right|^2dxd\tau \\ \notag
&&+\Int_{\bs 0}^t\Int_{\bs \Omega}\beta \left(-\Delta u_{\theta}(x,\tau)+q'(u_{\theta}(x,\tau))+(-\Delta)^{\frac{1}{2}}u_{\theta}(x,\tau)\right)^2 dxd\tau \\ \notag
&\leq&\Int_{\bs \Omega}\left(\frac{1}{2}|\nabla u_0(x)|^2+q(u_{0}(x))+u_0(x)(-\Delta)^{\frac{1}{2}}u_0\right) dx.
\end{eqnarray}
\end{proposition}

Proposition \ref{thm:ndg} is proved via  Galerkin approximations. Due to the presence of the stablizing function $g_{\theta}$, it is not obvious how to pass to the limit in the nonlinear term of the Galerkin approximations. Our main observation in this step is  strong convergence of  $\nabla u^N$ in $L^2(\Omega_T)$ which allows us to pass to the limit.

In order to obtain the weak solution to \eqref{eqn:u}, we consider the limit of $u_{\theta}$ as $\theta\rightarrow 0$.  The main difficulty is how to pass to the limit in the nonlinear term in the approximation equation. In \cite{DaiDu16},  the authors  proved the existence of weak solutions for degenerate Cahn-Hilliard equations by the following idea. The estimates for the positive mobility approximations yield  uniform bounds for $\partial_t u_{\theta_i}$ in $L^2(0,T; (H^2(\Omega))')$, and uniform bounds on $u_{\theta_i}$ in $L^{\infty}(0,T;H^1(\Omega)$.  Those uniform bounds yield  strong  convergence of $\sqrt{M_i(u_{\theta_i})}$ in $C(0,T;L^n(\Omega))$. By this and the weak convergence of $\sqrt{M_i(u_{\theta_i})}\nabla \mu_{\theta_i}$ in $L^2(\Omega_T)$, authors in  \cite{DaiDu16} showed (up to a subsequence) that  $M_{\theta_i}(u_{\theta_i})\nabla \mu_{\theta_i}\rightharpoonup \sqrt{M(u)}\xi $ weakly in $L^2(0,T;L^{\frac{2n}{n+2}}(\Omega))$ where $\xi$ is the weak limit of $\sqrt{M_i(u_{\theta_i})}\nabla \mu_{\theta_i}$.  The main task  left  is  to show $\sqrt{M(u)}\xi =M(u)(-\nabla \Delta u+q''(u)\nabla u)$ and the limit equation becomes a weak form Cahn-Hilliard equation.  Authors in \cite{DaiDu16}  proved that  this is almost true in the set where $u\neq \pm 1$. Their main idea is the following.  For small numbers $\delta_j$ monotonically decreasing to $0$, they  consider the limit in  a subset  $B_j$ of $\Omega_T$ where   approximate solutions converges uniformly and $|\Omega_T\backslash B_j|<\delta_j$.  By decomposing $B_j=D_j\cup\tilde D_j$ where mobility is bounded from below uniformly in $D_j$ and controlled above in $\tilde D_j$ by suitable multiples of $\delta_j$, they obtain the weak form equation for the limit function by passing  to the limit of $u_{\theta_i}$ on $D_j$ then letting $j$ goes to $\infty$.  Under further regularity assumptions on $\nabla \Delta u$, they obtained the explicit expression for $\xi$  in the weak form of the equation.

Due to the existence of the stablizing function $g(u)$ in our model,  it is much more delicate to carry out a similar analysis. The first obtacle is the bound estimate on $\partial_t u_{\theta_i}$ blows up when $\theta_i$ goes to zero and secondly,  it is more complicated to  derive an explicit expression of the weak limit of ${M_i(u_{\theta_i})}\nabla \frac{\mu_{\theta_i}}{g_{\theta_i}(u_{\theta_i})}$ in terms of $u$ in  the limit equation. We shall follow ideas in a recent work by the authors \cite{NXY22} by which  we  derive convergence of $g_{\theta_i}(u_{\theta_i})$ (consequently  $M_{\theta_i}(u_{\theta_i})$ from convergence of   $G_i(u)=\int_0^u g_{\theta_i}(s) ds$ . We then follow the idea in \cite{DaiDu16} to pass to the limit on the right hand side of the approximation equation. Below is our main theorem.

\begin{theorem}\label{thm:dg}
For any $u_0\in H^1(\Omega)$ and $T>0$, there exists a function $u:\Omega_T\rightarrow \R$ satisfying
\begin{itemize}
\item[i)] $u\in L^{\infty}(0,T;H^1(\Omega))\cap C([0,T];L^s(\Omega))\cap L^2(0,T; H^2(\Omega))$, where $1\leq s <\infty$,

\item[ii)] $g(u)\partial_tu\in L^{p}(0,T;(W^{1,q}(\Omega))')$ for $1\leq p<2$ and $q>2$.

\item[iii)] $u(x,0)=u_0(x)$ for all $x\in \Omega$,
\end{itemize}
which solves \eqref{eqn:u}-\eqref{eqn:mu} in the following weak sense
\begin{itemize}
\item[a)]There exists a set $B\subset \Omega_T$ with $|\Omega_T\backslash B|=0$ and a function $\zeta:\Omega_T\rightarrow \R^n$ satisfying $\chi_{B\cap P}M(u)\zeta\in L^{\frac{p}{p-1}}(0,T;L^{\frac{q}{q-1}}(\Omega,\R^n))$ such that
          \begin{eqnarray}\label{eqn:dg}
           &&\int_0^T<g(u)\partial_tu,\phi>_{(W^{1,q}(\Omega))',W^{1,q}(\Omega)}dt\\ \notag
&=&\int_{B\cap P}M(u)\zeta \cdot \nabla \phi dxdt-\int_{\Omega_T}\beta\left [\nabla  u\cdot \nabla \phi+q'(u)\phi+(-\Delta)^{\frac{1}{2}}(u)\phi\right]dxdt
          \end{eqnarray}  for all $\phi \in L^p(0,T;W^{1,q}(\Omega))$ with $p,q>2$. Here $P:=\{(x,t)\in \Omega_T: |1-u^2|\neq 0\}$ is the set where $M(u), g(u)$ are nondegenerate and $\chi_{B\cap P}$ is the characteristic function of set $B\cap P$.

\item[b)]  Assume $u\in L^2(0,T;H^2(\Omega)).$  For any open set $U\in \Omega_T$ on  which $g(u)>0 $ and $\nabla \Delta u \in L^p(U)$ for some $p>1$, we have
\begin{equation}
\zeta=\frac{-\nabla \Delta u+q''(u)\nabla u+\nabla (-\Delta)^{\frac{1}{2}}u}{g(u)}-\frac{g'(u)}{g^2(u)}\left(-\Delta u+q'(u)+(-\Delta)^{\frac{1}{2}} u \right)\nabla u . \label{eqn:zetaU}
\end{equation} a.e. in $U$.
\end{itemize}
Moreover, the following energy inequality holds for all $t>0$
\begin{eqnarray}
\label{eqn:ineq} &&\Int_{\bs \Omega}\left(\frac{1}{2}|\nabla u(x,t)|^2+q(u(x,t))+u(-\Delta)^{\frac{1}{2}}u\right)dz\\ \notag
&&+\Int_{\bs \Omega_r\cap B\cap P}M(u(x,\tau))|\zeta(x,\tau)|^2 dxd\tau \\ \notag
&&+\Int_{\bs \Omega_r\cap B\cap P}\beta\left(-\Delta u+q'(u)+(-\Delta)^{\frac{1}{2}}u\right)^2 dxd\tau \\ \notag
&\leq &\Int_{\bs \Omega}\left(\frac{1}{2}|\nabla u_0(x)|^2+q(u_0(x))\right)dx.
\end{eqnarray}
\end{theorem}

{  Lastly, we perform numerical simulations to validate our model. Using the proposed phase field model, we did simulations of  evolution of an elliptic prismatic loop and interactions between two circular prismatic loops under the combined effect of self-climb and non-conservative climb. Our numerial results indicate the self-climb effect slows down the shrinking of loop for the evolution of an elliptic prismatic loop. For interaction between two circular loops, the patterns in the two shrinking process are quite different with or without  the self-climb effect .}

The paper is organized as follows. We shall derive sharp interface limit for \eqref{eqn:ch} and \eqref{eqn:chem} through formal asymptotic expansions in section \ref{sec:asym}. Section \ref{sec:ndg} is devoted to the proof of Theorem \ref{thm:ndg} and Theorem\ref{thm:dg}  is proved in section \ref{sec:dg}. {  Numerical simulations are presented in section \ref{sec:simulation}}.

\section{Sharp interface limit via asymptotic expansions}
\label{sec:asym}
In this section, we perform a formal asymptotic analysis to obtain the dislocation self-climb velocity of the proposed phase field model \eqref{eqn:ch} and \eqref{eqn:chem} in the sharp interface limit $\e \rightarrow 0$.
\subsection{Outer expansions}
We first perform expansion in the region far from the dislocations. Assume the expansion for $u$ is
\begin{equation}
u(x,y,t)=u^{(0)}(x,y,t)+u^{(1)}(x,y,t)\e+u^{(2)}(x,y,t)\e^2+\hdots
\end{equation}
Correspondingly, we have
\begin{eqnarray*}
&&M(u)=M(u^{(0)})+M'(u^{(0)})u^{(1)}\e +\left(M'(u^{(0)})u^{(2)}+\frac{1}{2}M''\left (u^{(0)}\right)\left(u^{(1)}\right)^2\right)\e^2+\hdots,\\
&&g(u)=g(u^{(0)})+g'(u^{(0)})u^{(1)}\e +\left(g'(u^{(0)})u^{(2)}+\frac{1}{2}g''(u^{(0)})\left(u^{(1)}\right)^2\right)\e^2+\hdots,\\
&&q'(u)=q'(u^{(0)})+q''(u^{(0)})u^{(1)}\e +\left(q''(u^{(0)})u^{(2)}+\frac{1}{2}q^{(3)}(u^{(0)})\left(u^{(1)}\right)^2\right)\e^2+\hdots,\\
&&f^d_{cl}(x,y,u)=f^d_{cl}(x,y,u^{(0)})+f^d_{cl}(x,y,u^{(1)})\e+f^d_{cl}(x,y,u^{(2)})\e^2+\hdots.
\end{eqnarray*}
We also expand the chemical potential $\mu$ as
\begin{equation}
\mu=\frac{1}{\e^2}\left(\mu^{(0)}+\mu^{(1)}\e+\mu^{(2)}\e^2+\hdots\right). \label{eqn:muexp}
\end{equation}
Rewrite equation \eqref{eqn:ch} as
\begin{equation}
g(u)(\partial_t u+\beta\mu)=M_0\nabla \cdot (\nabla \mu-\mu \frac{g'(u)}{g(u)}\nabla u), \label{eqn:ch1}
\end{equation}
and set
\begin{equation}
w=-\mu \frac{g'(u)}{g(u)}=\frac{1}{\e^2}\left(w^{(0)}+w^{(1)}\e+w^{(2)}\e^2+\hdots\right).
\end{equation}

Plugging the expansions into \eqref{eqn:ch1} and \eqref{eqn:chem} and matching the coefficients of $\e$ powers in both equations, the $O(\frac{1}{\e^2})$ equations of \eqref{eqn:ch1}  and \eqref{eqn:chem} yield
\begin{eqnarray}
\label{eqn:u0}\beta g(u^{(0)})\mu^{(0)}&=&M_0\nabla \cdot  \left(\nabla \mu^{(0)}+w^{(0)}\nabla u^{(0)}\right),\\
\label{eqn:mu0}\mu^{(0)}&=&q'(u^{(0)}).
\end{eqnarray}
Since
$$
 w^{(0)}=\mu^{(0)}\frac{ g'(u^{(0)})}{g(u^{(0)})},
$$
then $u^{(0)}=1$ or $u^{(0)}=-1$ satisfies equations \eqref{eqn:u0}-\eqref{eqn:mu0}. In particular, such choice of $u^{(0)}$ implies $\mu^{(0)}=0$.

The $O(\frac{1}{\e})$ equations of \eqref{eqn:ch1}  and \eqref{eqn:chem} yield
\begin{eqnarray}
\label{eqn:u1}\beta \left( g(u^{(0)})\mu^{(1)}+g'(u^{(0)})u^{(1)}\mu^{(0)}\right)&=&M_0\nabla \cdot  \left(\nabla \mu^{(1)}+w^{(0)}\nabla u^{(1)}+w^{(1)}\nabla u^{(0)}\right),\\
\label{eqn:mu1}\mu^{(1)}&=&q''(u^{(0)})u^{(1)}+h(u^{(0)})f^d_{cl}(x,y,u^{(0)}).
\end{eqnarray}
Since $u^{(0)}=1$ or $-1$,
 $u^{(1)}=0$ satisfies \eqref{eqn:u1}-\eqref{eqn:mu1}. Moreover, such choice of $u^{(1)}$ guarantees $\mu^{(1)}=0$.

The $O(1)$ equations of \eqref{eqn:ch1}  and \eqref{eqn:chem} are
\begin{eqnarray*}
&&u_t^{(0)}g(u^{(0)})+\beta \left(g(u^{(0)})\mu^{(2)}+g'(u^{(0)})u^{(1)}\mu^{(1)}+\mu^{(0)}\left(g'(u^{(0)})u^{(2)}+\frac{1}{2}g''(u^{(0)})\left(u^{(1)}\right)^2\right)\right)\\
&=&M_0\nabla \cdot  \left(\nabla \mu^{(2)}+w^{(0)}\nabla u^{(2)}+w^{(1)}\nabla u^{(1)}+w^{(2)}\nabla u^{(0)}\right),\\
\mu^{(2)}&=&-\Delta u^{(0)}+q''(u^{(0)})u^{(2)}+\frac{1}{2}q^{(3)}(u^{(0)})(u^{(1)})^2\\
 &&+h(u^{(0)})f^d_{cl}(x,y,u^{(1)})+h'(u^{(0)})u^{(1)}f^d_{cl}(x,y,u^{(0)}).
\end{eqnarray*}
 Taking into account of the fact $u^{(0)}=\pm 1$, $u^{(1)}=\mu^{(0)}=\mu^{(1)}=0$, the equations above reduce to
\begin{eqnarray}
\label{eqn:u2}0
&=&\nabla \cdot  \left(\nabla \mu^{(2)}+w^{(0)}\nabla u^{(2)}\right),\\
\label{eqn:mu2}\mu^{(2)}&=&q''(u^{(0)})u^{(2)}
+h(u^{(0)})f^d_{cl}(x,y,u^{(1)}).
\end{eqnarray}
Thus $u^{(2)}=0$ satisfies  \eqref{eqn:u2}-\eqref{eqn:mu2}. Moreover, such choice of $u^{(2)}$ guarantees $\mu^{(2)}=0$.

In general, if $u^{(0)}=\pm 1$, $u^{(1)}=u^{(2)}=\hdots=u^{(k+1)}=0$, the $O(\e^k)$ of the $k\geq 1$ equations of \eqref{eqn:ch1}  and \eqref{eqn:chem} yield
\begin{eqnarray}
\label{eqn:uk}  0&=&\nabla \cdot  \left(\nabla \mu^{(k+2)}+w^{(0)}\nabla u^{(k+2)}\right),\\
\label{eqn:muk}  \mu^{(k+2)}&=&q''(u^{(0)})u^{(k+2)}+h(u^{(0)})f^d_{cl}(x,y,u^{(k+1)}).
\end{eqnarray}
Thus $u^{(k+2)}=0$ satisfies \eqref{eqn:uk} and \eqref{eqn:muk}.

In summary, we have $u=1$ or $u=-1$ in the outer region.

\subsection{Inner expansions}
For the small inner regions near the dislocations, we introduce local coordinates near the dislocations. Considering a dislocation $C$ parameterized by its arc length parameter $s$. We denote a point on the dislocation by $\rr_0(s)$ with  tangent unit vector $\mathbf{t}(s)$ and inward normal vector $\bn(s)$. A point near the dislocation $C$ is expressed as
\begin{equation}
\rr(s,d)=\rr_0(s)+d\mathbf{n}(s),
\end{equation}
where $d$ is the signed distance from point $\rr$ to the dislocation. Since the gradient fields are of order $O(\frac{1}{\e})$, we introduce the variable $\rho=\frac{d}{\e}$ and use coordinates $(s,\rho)$ in the inner region.
Under this setting, we write $u(x,y,t)=U(s,\rho,t)$ and equation \eqref{eqn:ch} can be written as
\begin{eqnarray}
&&\label {eqn:chin}\hspace{0.5 in} g(U)\left(\partial_t U-\frac{1}{\e}v_n\partial_{\rho}U+\beta \mu\right)=\frac{M_0}{1-\e\rho\kappa}\partial_s\left(\frac{1}{1-\e\rho\kappa}\left(\partial_s \mu-\mu\frac{g'(U)}{g(U)}\partial_s U\right)\right) \\
&& \notag \hspace{1.9 in}+\frac{1}{\e^2}\frac{M_0}{1-\e\rho\kappa}\partial_{\rho}\left((1-\e\rho\kappa)\left(\partial_{\rho}\mu-\mu \frac{g'(U)}{g(U)}\partial_{\rho}U\right)\right),\\
&&\label{eqn:chemin} \hspace{1 in} \mu=-\frac{1}{1-\e\rho\kappa}\partial_s\left(\frac{1}{1-\e\rho\kappa}\partial_sU\right)-\frac{1}{\e^2}\frac{1}{1-\e\rho\kappa}\partial_{\rho}\left((1-\e\rho\kappa)\partial_{\rho}U\right)\\
&& \notag \hspace{1.2in}+\frac{1}{\e^2}q'(U)+\frac{1}{\e}h(U)f_{cl}(s,\rho,U).
\end{eqnarray}

Assume that $\mu$ takes the same form expansion as  \eqref{eqn:muexp}. The following expansions hold for $U$ and the climb force $f_{cl}$ within dislocation core region:
\begin{equation}
U(s,\rho,t)=U^{(0)}(\rho)+\e U^{(1)}(s,\rho,t)+\e^2 U^{(2)}(s,\rho,t)+\hdots, \label{eqn:U}
\end{equation}
and
\begin{equation}
f_{cl}(s,\rho,U)=\frac{1}{\e}f^{(-1)}_{cl}(\rho,U)+f^{(0)}_{cl}(s)+O(\e), \label{eqn:fcl}
\end{equation}
where
\begin{eqnarray}
 f^{(-1)}_{cl}(\rho,U)&=&\frac{G b^2}{4\pi(1-\nu)}\Int_{\bs-\infty}^{+\infty}\frac{\partial_{\rho}U(\tau)}{\rho-\tau}d\tau,\\
 f^{(0)}_{cl}(s)&=&f^d_{cl}(s)+f^{app}_{cl}(s),\\ \label{eqn:fdcl}
 f^d_{cl}(s)&=&\frac{G b^2}{4\pi(1-\nu)}\kappa \ln \e+O(1).
\end{eqnarray}
Here we assume the leading order solution $U^{(0)}$, which describe the dislocation core profile, remains the same at all points on the dislocation at any time. The term $ \frac{1}{\e}f^{(-1)}_{cl}(\rho,U)$ in the climb force expansion is due to the singular stress field near the dislocation and vanishes on the dislocation (i.e. $f^{(-1)}_{cl}(\rho,U^{(0)})=0$). The climb force $f^d_{cl}(s)$ is generated by the dislocations and has asymptotic expansions \eqref{eqn:fdcl}. This asymptotic expansion of climb force $f_{cl}$ in the phase field model was obtained in \cite{Niu2021} based on dislocation theories \cite{Hirth-Lothe,GB,ZhaoXiang}.

Letting \begin{equation} W=\mu\frac{g'(U)}{g(U)}=\frac{1}{\e^2}\left(W^{(0)}+W^{(1)}\e+W^{(2)}\e^2+\hdots \right),\end{equation}
the leading orders of equations \eqref{eqn:chin} and \eqref{eqn:chemin} are $O(\frac{1}{\e^4})$ and $O(\frac{1}{\e^2})$, respectively, which yield
\begin{eqnarray}
&&\label{eqn:U0} 0= \partial_{\rho}\left(\partial_{\rho}\mu^{(0)}-W^{(0)}\partial_{\rho}U^{(0)}\right),\\
&& \label{eqn:mu0in} \mu^{(0)}=-\partial_{\rho\rho}U^{(0)}+q'(U^{(0)})+h(U^{(0)})f^{(-1)}_{cl}(\rho,U^{(0)}).
\end{eqnarray}
Integrating Eq.~\eqref{eqn:U0}, we have
\begin{equation}
\partial_{\rho}\mu^{(0)}-W^{(0)}\partial_{\rho}U^{(0)}=C_0(s). \label{eqn:c0}
\end{equation}
Since $\mu^{(0)}=0$,  $u^{(0)}=1$ or $-1$ in the outer region, we must have $\mu^{(0)}\rightarrow 0$ and $\partial_{\rho}U^{(0)}\rightarrow 0$ as $\rho\rightarrow \pm \infty$. Therefore $C_0(s)=0$.   Dividing \eqref{eqn:c0} by $\mu^{(0)}$ and taking integration, using $W^{(0)}=\mu^{(0)}\frac{g'(U^{(0)})}{g(U^{(0)})}$, we have $\mu^{(0)}=\tilde C_0(s)g(U^{(0)})$.
Since $\mu^{(0)}/g(u^{(0)})$ is $\infty$ in the outer region,
we must have $\tilde C_0(s)=0$.
Thus
\begin{equation}
\mu^{(0)}=-\partial_{\rho\rho}U^{(0)}+q'(U^{(0)})+h(U^{(0)})f^{(-1)}_{cl}(\rho, U^{(0)})=0. \label{eqn:Umu0}
\end{equation}
Solution $U^{(0)}$  to \eqref{eqn:Umu0} subject to far field condition $U^{(0)}{(+\infty)}=-1$ and $U^{(0)}{(-\infty)}=1$ can be found numerically (see \cite{Niu2021}). In particular, $\partial_{\rho}U^{(0)}<0$ for all $\rho$.

Next, the $O(\frac{1}{\e^3})$ equation of \eqref{eqn:chin} and $O(\frac{1}{\e})$ equation of  \eqref{eqn:chemin} yield, using $\mu^{(0)}=0$, that
\begin{eqnarray}
\label{eqn:U1} 0&=& \partial_{\rho}\left(\partial_{\rho}\mu^{(1)}-W^{(1)}\partial_{\rho}U^{(0)}\right),\\
 \label{eqn:mu1in} \mu^{(1)}&=&-\partial_{\rho\rho}U^{(1)}+\kappa \partial_{\rho}U^{(0)}+q''(U^{(0)})U^{(1)}+h'(U^{(0)})f^{(-1)}(\rho, U^{(0)})U^{(1)}\\ \notag &&
 +h(U^{(0)})f^{(0)}_{cl}(s).
\end{eqnarray}
Similar to the calculation from Eq.~\eqref{eqn:U0} to Eq.~\eqref{eqn:c0} given above, by matching with the outer solutions, we have $\partial_{\rho}\mu^{(1)}-W^{(1)}\partial_{\rho}U^{(0)}=0$.
When $\mu^{(0)}=0$, we have  $W^{(1)}=\mu^{(1)}\frac{g'(U^{(0)})}{g(U^{(0)})}$,
the obtained equation becomes
\begin{equation}
\partial_{\rho}\mu^{(1)}-\mu^{(1)}\partial_{\rho}\ln g(U^{(0)})=0. \label{eqn:c1}
\end{equation}
Dividing  \eqref{eqn:c1} by $\mu^{(1)}$ and integrating, we have $\mu^{(1)}=\tilde C_1(s)g(U^{(0)})$. Thus equation \eqref{eqn:mu1in} can be rewritten as
\begin{equation}
LU^{(1)}=-\kappa \partial_{\rho}U^{(0)}-h(U^{(0)})f^{(0)}_{cl}(s)+\tilde C_1(s)g(U^{(0)}),\label{eqn:L}
\end{equation}
where $L=-\partial_{\rho\rho}+q''(U^{(0)})+h'(U^{(0)})f^{(-1)}(\rho, U^{(0)})$. 
Multiplying both sides of Eq.~\eqref{eqn:L} by $\partial_{\rho}U^{(0)}$ and integrate with respect to $\rho$ over $(-\infty, +\infty)$, we have
$$
\int_{-\infty}^{+\infty}\left(-\kappa \partial_{\rho}U^{(0)}-h(U^{(0)})f^{(0)}_{cl}(s)+\tilde C_1(s)g(U^{(0)})\right)\partial_{\rho}U^{(0)} d\rho=0.
$$
From this, we conclude
$$
\tilde C_1(s)=-\alpha \kappa+H_0 f^{(0)}_{cl}(s),
$$
where   $\alpha>0$ is  given by
$$
\alpha=-\frac{\int_{-\infty}^{+\infty}\left(\partial_{\rho}U^{(0)}\right)^2 d\rho}{\int_{-\infty}^{+\infty}g(U^{(0)})\partial_{\rho}U^{(0)}d\rho}.
$$
Therefore
\begin{equation}
\mu^{(1)}=g(U^{(0)})\left(-\alpha \kappa+H_0 f^{(0)}_{cl}(s)\right). \label{eqn:mu1exp}
\end{equation}

Letting $\overline{\mu}=\frac{\mu}{g(U)}$, \eqref{eqn:chin} can be written as
\begin{eqnarray}\notag
&& g(U)\left(\partial_t U-\frac{1}{\e}v_n\partial_{\rho}U+\beta \mu\right)\\ \label{eqn:chinmubar}
&=& \frac{M_0}{1-\e\rho\kappa}\partial_s\left(\frac{g(U)}{1-\e\rho\kappa}\left(\partial_s \overline \mu\right) \right)+\frac{1}{\e^2}\frac{M_0}{1-\e\rho\kappa}\partial_{\rho}\left(\left(1-\e\rho\kappa\right)g(U)\partial_{\rho}\overline\mu\right)
\end{eqnarray}
Using $\mu^{(0)}=0$, $\partial_{\rho}\overline \mu^{(1)}=\partial_{\rho}\frac{\mu^{(1)}}{g(U^{(0)})}=0$, the $O(\frac{1}{\e^2})$ order equation of \eqref{eqn:chinmubar} reduces to
$$
\partial_{\rho}\left(g(U^{(0)})\partial_{\rho}\overline{\mu}^{(2)}\right)=0.
$$
Integrating with respect to $\rho$, we have $g(U^{(0)})\partial_{\rho}\overline{\mu}^{(2)}=C_2(s)$. Matching with outer solutions, we must have $C_2(s)=0$. Thus $\partial_{\rho}\overline{\mu}^{(2)}=0$ which gives  $\overline{\mu}^{(2)}=\tilde C_2(s)$.

Next we look at the $O(\frac{1}{\e})$ equation of \eqref{eqn:chinmubar}. Using $\mu^{(0)}=0$, $\partial_{\rho}\overline \mu^{(1)}=0$ and $\partial_{\rho}\overline \mu^{(2)}=0$, we have
$$
g(U^{(0)})(-v_n\partial_{\rho} U^{(0)}+\beta \mu^{(1)})=M_0\partial_s\left(g(U^{(0)})\partial_s\overline \mu^{(1)}\right)+M_0\partial_{\rho}\left(g(U^{(0)})\partial_{\rho} \overline\mu^{(3)}\right).
$$
Integrating this equation with respect to $\rho$ and matching with outer solutions yields
\begin{equation}
v_n=\lambda \partial_{ss}\overline \mu^{(1)} -\eta \overline\mu^{(1)}\label{eqn:vn}
\end{equation}
where we used the fact that $g(U^{(0)})$ is independent of $s$, $\overline\mu^{(1)}=-\alpha \kappa+H_0 f^{(0)}_{cl}(s)$ by \eqref{eqn:mu1exp},
 and
\begin{equation}\label{eqn:lambda-eta}
\lambda=-\frac{M_0\int_{-\infty}^{+\infty}g(U^{(0)})d\rho}{\int_{-\infty}^{+\infty}g(U^{(0)})\partial_{\rho}U^{(0)}d\rho}>0, \hspace{0.1in} \eta=-\frac{\beta\int_{-\infty}^{+\infty}g(U^{(0)})d\rho}{\int_{-\infty}^{+\infty}g(U^{(0)})\partial_{\rho}U^{(0)}d\rho}>0.
\end{equation}
Substitute $\overline\mu^{(1)}=-\alpha \kappa+H_0 f^{(0)}_{cl}(s)$
into \eqref{eqn:vn},  the sharp interface limit equation is
\begin{equation}
v_n=-\lambda\partial_{ss}\left(\alpha \kappa-H_0 f^{(0)}_{cl}(s)\right)+\eta \left(\alpha \kappa-H_0 f^{(0)}_{cl}(s)\right). \label{eqn:sharpif}
\end{equation}

\begin{remark}
The velocity in the obtained sharp interface limit equation
\eqref{eqn:sharpif}  is a combination of the dislocation self-climb velocity \cite{Niu2017,Niu2019,Niu2021} (the first term), and the dislocation climb velocity by mobility law \cite{Xiang2003,Xiang2006,Arsenlis2007} (the second term). The coefficients of these two contributions are determined through Eq.~\eqref{eqn:lambda-eta} by the parameters $M_0$ and $\beta$, respectively, in the phase field model in \eqref{eqn:ch} based on the physics. Note that the curvature term in both contributions is a correction to the dislocation self-force to fix the problem of larger numerical dislocation core size in the phase field model than the actual dislocation core size \cite{Niu2021}. We have mentioned previously that the factor $g(U)$ on the left-hand side in Eq.~\eqref{eqn:model0} is mainly for the wellposedness proof, and without it, the dislocation velocity given by the sharp interface limit is similar, with $\lambda=\frac{M_0}{2}\int_{-\infty}^{+\infty}g(U^{(0)})d\rho$ and $\eta=\frac{\beta}{2}\int_{-\infty}^{+\infty}g(U^{(0)})d\rho$.
\end{remark}

\section {Weak solution for phase field model}
\subsection {Weak solution for phase field model with positive mobilities}
\label{sec:ndg}
In this subsection, we prove existence of weak solutions for phase field model with positive mobilities summarized in Proposition~\ref{thm:ndg}.

Let $\mathbb Z_{+}$ be the set of nonnegative integers and $\Omega=[0,2\pi]^n$ with $n\leq 2$.  We choose an orthonormal basis for $L^2(\Omega)$ as
\begin{eqnarray*}
\{\phi_j:j=1,2,\hdots\}=\left \{(2\pi)^{-{n/2}},\text{Re}\left(\pi^{-n/2}e^{i\xi\cdot x}\right), \text{Im}\left(\pi^{-n/2}e^{i\xi\cdot x}\right):\xi \in \mathbb Z^n_{+}\backslash \{0,\hdots,0\}\right \}.
\end{eqnarray*}
Observe $\{\phi_j\}$ is also orthogonal in $H^k(\Omega)$ for any $k\geq 1$.

\subsubsection{Galerkin approximations}

Define
$$
u^N(x,t)=\sum_{j=1}^Nc_j^N(t)\phi_j(x), \hspace{0.5 in} \mu^N(x,t)=\sum_{j=1}^Nd^N_{j}(t)\phi_j(x),
$$
where $\{c_j^N,d^N_j\}$ satisfy
\begin{eqnarray}\label{eqn:uN}
\Int_{\bs \Omega}\partial_t u^N\phi_j dx&=&-\Int_{\bs \Omega}M_{\theta}(u^N)\nabla \frac{\mu^N}{g_{\theta}(u^N)}\cdot \nabla \frac{\phi_j}{g_{\theta}(u^N)}dx-\beta\int_{\Omega}\mu^N \phi_j dx\\ \label{eqn:muN}
\Int_{\bs \Omega}\mu^N\phi_j dx&=&\Int_{\bs \Omega}\left(\nabla u^N\cdot \nabla \phi_j+q'(u^N)\phi_j+\phi_j(-\Delta)^{\frac{1}{2}} u^N \right)dx,\\ \label{eqn:uN0}
u^N(x,0)&=&\sum_{j=1}^N\left(\int_{\Omega} u_0\phi_j dx\right)\phi_j(x).
\end{eqnarray}
\eqref{eqn:uN}-\eqref{eqn:uN0} is an initial value problem for a system of ordinary equations for $\{c_j^N(t)\}$. Since right hand side of \eqref{eqn:uN} is continuous in $c_j^N$, the system has a local solution.

Define energy functional
$$
E(u)=\Int_{\bs \Omega}\left\{\frac{1}{2}|\nabla u|^2 +q(u)+|(-\Delta)^{\frac{1}{4}}u|^2\right\} dx.
$$
Direct calculation using \eqref{eqn:uN} and \eqref{eqn:muN} yields
$$
\frac{d}{dt}E(u^N(x,t))=-\scaleint{7ex}_{\bs\Omega}M_{\theta}(u^N)\left|\nabla \frac{\mu^N}{g_{\theta}(u^N)}\right|^2 dx-\beta\int_{\Omega}\left(\mu^N\right)^2 dx,
$$
integration over $t$ gives the following energy identity.
\begin{eqnarray}\notag
&&\scaleint{7ex}_{\bs \Omega}\left(\frac{1}{2}|\nabla u^N(x,t)|^2+q(u^N(x,t))+u^N(-\Delta)^{\frac{1}{2}}u^N\right)dx\\ \notag
&&+\Int_{\bs 0}^t\Int_{\bs\Omega}\left[M_{\theta}(u^N(x,\tau))\left|\nabla \frac{\mu^N(x,\tau)}{g_{\theta}(u^N(x,\tau))}\right|^2 +\beta \left(\mu^N\right)^2\right]dxd\tau\\ \label{eqn:energyid}
&=&\Int_{\bs \Omega}\left(\frac{1}{2}|\nabla u^N(x,0)|^2+q(u^N(x,0))+u^N(x,0)(-\Delta)^{\frac{1}{2}}u^N(x,0)\right)dx\\ \notag
&\leq&\left \Arrowvert \nabla u_0\right \Arrowvert_{L^2(\Omega)}^2+C\left(\left\Arrowvert u_0\right\Arrowvert^{r+1}_{H^1{\Omega}}+|\Omega|\right)+\frac{1}{2}\left\Arrowvert u_0\right\Arrowvert_{L^2(\Omega)}^2\leq C <\infty
\end{eqnarray}
Here and throughout the paper, $C$ represents a generic constant possibly depending only on  $\beta$,  $T$, $\Omega$,  $u_0$ but not on $\theta$. Since $\Omega$ is bounded region, by  growth assumption  \eqref{eqn:qu} and Poincare's inequality, the energy identity \eqref{eqn:energyid} implies
$
u^N\in L^{\infty}(0,T;H^1(\Omega))
$
and $\mu^N\in L^2(\Omega_T)$ with
\begin{equation}
\left\Arrowvert \mu^N\right\Arrowvert _{L^2(\Omega_T)}, \left\Arrowvert u^N\right\Arrowvert _{L^{\infty}(0,T;H^1(\Omega))}\leq C \text{ for all } N, \label{eqn:uNbd}
\end{equation}
and
\begin{equation}
\left \Arrowvert \sqrt{M_{\theta}(u^N)}\nabla \frac{\mu^N}{g_{\theta}(u^N)}\right \Arrowvert_{L^2(\Omega_T)}\leq C \text{ for all } N.  \label{eqn:muNbd}
\end{equation}
By \eqref{eqn:uNbd}, the coefficients $\{c_j^N(t)\}$ are bounded in time, thus the system \eqref{eqn:uN}-\eqref{eqn:uN0} has a global solution. In addition, by Sobolev embedding theorem and growth assumption \eqref{eqn:q'u} on $q'(u)$, we have
$$
q'(u^N)\in L^{\infty}(0,T; L^p(\Omega)), \hspace{0.5in} M_{\theta}(u^N)\in L^{\infty}(0,T;L^p(\Omega))
$$
 for any $1\leq p<\infty$ with
\begin{eqnarray} \label{eqn:qunbd}
&&\left \Arrowvert q'(u^N)\right\Arrowvert_{L^{\infty}(0,T;L^p(\Omega))}\leq C  \text{ for all } N, \\ \label{eqn:mbd}
&&\left \Arrowvert M_{\theta}(u^N)\right \Arrowvert_{L^{\infty}(0,T;L^p(\Omega))}\leq C     \text{ for all } N.
\end{eqnarray}
\subsubsection{Convergence of $u^N$}\label{sec:uNconv}

Given $q>2$ and any $\phi\in L^2(0,T;W^{1,q}(\Omega))$, let  $\Pi_N\phi(x,t)=\sum_{j=1}^N\left(\int_{\Omega}\phi(x,t)\phi_j(x)dx\right) \phi_j(x)$ be the orthogonal projection of $\phi$ onto span$\{\phi_j\}_{j=1}^N$.  Then
\begin{eqnarray*}
&&\left \arrowvert \int_{\Omega}\partial_tu^N\phi dx\right \arrowvert=\left\arrowvert \int_{\Omega}\partial_tu^N\Pi_N\phi  dx\right \arrowvert\\
&=&\left \arrowvert \Int_{\bs \Omega}\left[M_{\theta}(u^N)\nabla \frac{\mu^N}{g_{\theta}(u^N)}\cdot \nabla \frac{\Pi_N\phi}{g_{\theta}(u^N)}-\beta \mu^N\Pi_N\phi \right]dx\right \arrowvert\\
&\leq& \left\Arrowvert \sqrt{M_{\theta}(u^N)} \nabla \frac{\mu^N}{g_{\theta}(u^N)}\right \Arrowvert_{L^2(\Omega)} \left\Arrowvert \sqrt{M_{\theta}(u^N)} \nabla\frac{\Pi_N\phi}{g_{\theta}(u^N)}\right \Arrowvert_{L^2(\Omega)}+\beta||\mu^N||_{L^2(\Omega)}||\phi||_{L^2(\Omega)}.
\end{eqnarray*}

Since
$$
 \nabla\frac{\Pi_N\phi}{g_{\theta}(u^N)}=\frac{1}{g_{\theta}(u^N)}\nabla \Pi_N \phi-\Pi_N \phi \frac{g_{\theta}'(u^N)}{g^2_{\theta}(u^N)}\nabla u^N,
$$
we have
\begin{eqnarray*}
&& \Int_{\bs \Omega}M_{\theta}(u^N)\left \arrowvert \nabla\frac{\Pi_N\phi}{g_{\theta}(u^N)}\right \arrowvert^2 dx\\
&\leq& 2 M_0\Int_{\bs\Omega}\left(\frac{1}{g_{\theta}(u^N)}\left \arrowvert \nabla \Pi_N \phi\right \arrowvert ^2 + \frac{|g'_{\theta}(u^N)|^2}{g^3_{\theta}(u^N)} | \Pi_N \phi|^2|\nabla u^N|^2\right)dx\\
&\leq& C(M_0,\theta)\left(\left \Arrowvert \nabla \Pi_N \phi \right \Arrowvert^2_{L^2(\Omega)}+ \left \Arrowvert \Pi_N \phi \right\Arrowvert^2_{L^{\infty}(\Omega)}\left \Arrowvert \nabla u^N \right \Arrowvert^2_{L^2(\Omega)}\right)\\
&\leq& C(M_0,\theta)\left(\left \Arrowvert \Pi_N \phi \right \Arrowvert^2_{W^{1,q}(\Omega)}\right) \leq C(M_0,\theta)\left \Arrowvert\phi \right \Arrowvert^2_{W^{1,q}(\Omega)}.
\end{eqnarray*}
Therefore
\begin{equation}
\left \Arrowvert \partial_t u^N \right \Arrowvert_{L^2(0,T; (W^{1,q}(\Omega))')} \leq C(M_0,\theta) \text{ for all } N. \label{eqn:uNtbd}
\end{equation}

 For $1\leq s < \infty$, since $n\leq 2$, by Sobolev embedding theorem and Aubin-Lions Lemma (see \cite{Sim86} and  Remark \ref{rmk:cpt}) , the following embeddings are compact :
$$
\left \{f\in L^2(0,T;H^1(\Omega)):\partial_tf\in L^2(0,T;(W^{1,q}(\Omega))'\right \}\hookrightarrow  L^2(0,T;L^s(\Omega)),
$$
and
$$
\left \{f\in L^{\infty}(0,T;H^1(\Omega)):\partial_tf\in L^2(0,T;(W^{1,q}(\Omega))'\right \}\hookrightarrow  C([0,T];L^s(\Omega)).
$$
From this  and  the boundedness of $\{u^N\}$ and $\{\partial_t u^N\}$, we can find a subsequence, and $u_{\theta} \in L^{\infty}(0,T;H^1(\Omega))$ such that as $N\rightarrow \infty$, for $1\leq s< \infty$.
\begin{eqnarray}\label{eqn:uNwC}
u^N &\rightharpoonup& u_{\theta}\text{  weakly-* in }L^{\infty}(0,T;H^1(\Omega)), \\ \label{eqn:uNC}
u^N &\rightarrow& u_{\theta} \text{ strongly in } C([0,T];L^s(\Omega)),\\ \label{eqn:uNp}
u^N&\rightarrow &
u_{\theta} \text{ strongly in } L^2(0,T;L^s(\Omega)) \text{ and  a.e.  in }\Omega_T, \\\label{eqn:uNtC}
\partial_tu^N &\rightharpoonup &\partial_t u_{\theta} \text{ weakly in } L^2(0,T;(W^{1,q}(\Omega))').
\end{eqnarray}
 In addition
$$
\left \Arrowvert u_{\theta} \right \Arrowvert_{L^{\infty}(0,T;H^1(\Omega))} \leq C, \hspace{0.5 in} \left \Arrowvert \partial_t u_{\theta}\right \Arrowvert_{L^2(0,T; (W^{1,q}(\Omega))')} \leq C(M_0, \theta).
$$

By \eqref{eqn:uNC}, growth assumption \eqref{eqn:q'u} on $q'(u^N)$, and general dominated convergence Theorem, we have
\begin{eqnarray}\label{eqn:MuNC}
&& M_{\theta}(u^N) \rightarrow M_{\theta}(u_{\theta}) \text{ strongly in } C([0,T]; L^s(\Omega))\\ \label{eqn:sqMC}
&& \sqrt {M_{\theta}(u^N)} \rightarrow \sqrt {M_{\theta}(u_{\theta})} \text{ strongly in } C([0,T]; L^s(\Omega))\\ \label{eqn:quNC}
&& q'(u^N) \rightarrow q'(u_{\theta}) \text{ strongly in } C([0,T];L^s(\Omega))
\end{eqnarray}
for $1\leq s<\infty$. By \eqref{eqn:qunbd} and \eqref{eqn:quNC}, we have
\begin{equation}
q'(u^N) \rightharpoonup q'(u_{\theta}) \text{ weakly-* in } L^{\infty}([0,T];L^s(\Omega)). \label{eqn:quNwC}
\end{equation}
\begin{remark}\label{rmk:cpt}
Let $X$, $Y$, $Z$ be Banach spaces with compact embedding $X\hookrightarrow Y$ and continuous embedding $Y\hookrightarrow Z$. Then the embeddings
\begin{equation}
\{f \in L^p(0,T;X); \partial_t f \in L^1(0,T;Z)\} \hookrightarrow L^p(0,T;Y) \label{eqn:l2embed}
\end{equation}
and
\begin{equation}
\{f \in L^{\infty}(0,T;X); \partial_t f \in L^r(0,T;Z)\} \hookrightarrow C([0,T];Y) \label{eqn:linftyembed}
\end{equation}
are compact for any $1\leq p < \infty$ and $r>1$ (Corollary 4, \cite{ Sim86}, see also \cite{Lions69}) . For convergence of $u^N$, we apply this for $p=2=r$ with  $X=H^1(\Omega)$, $Y=L^s(\Omega)$ for $1\leq s<\infty$ and $Z=W^{1,q}(\Omega)'$.
\end{remark}
\subsubsection{Weak solution}
By \eqref{eqn:muN}, we have
$$
\int_{\Omega} \mu^N u^N dx=\int_{\Omega}\left(|\nabla u^N|^2 dx +q'(u^N)u^N+u^N (-\Delta)^{\frac{1}{2}}u^N\right)dx.
$$
Integration with respect to $t$ from $0$ to $T$  gives
\begin{eqnarray*}
&&\int_{\Omega_T} \mu^N(x,\tau) u^N(x,\tau) dxd\tau \\
&=&\int_{\Omega_T}\left(\nabla u^N(x,\tau)|^2 dx +q'(u^N(x,\tau))u^N(x,\tau)+u^N (-\Delta)^{\frac{1}{2}}u^N\right)dxd\tau.
\end{eqnarray*}

By \eqref{eqn:uNbd}. there exists a subsequence of $\mu^N$, not relabeled, converges weakly to $\mu_{\theta} \in L^2(\Omega_T)$.
Passing to the limit in the equation above, by \eqref{eqn:uNp}, \eqref{eqn:quNC},  we have
\begin{eqnarray} \label{eqn:guNC}
\int_{\Omega_T} \mu_{\theta}u_{\theta} dxd\tau&=&\lim_{N\rightarrow \infty}\int_{\Omega_T}|\nabla u^N|^2 dxd\tau+\int_{\Omega_T}q'(u_{\theta})u_{\theta} dxd\tau\\ \notag
&&+\int_{\Omega_T}u_{\theta}(-\Delta)^{\frac{1}{2}}u_\theta dx d\tau
\end{eqnarray}
On the other hand,
\begin{eqnarray}\label{eqn:muNuN}
&&\hspace {0.2 in} \int_{\Omega_T}\mu^N(x,\tau) u_{\theta}(x,\tau)dxd\tau=\int_{\Omega_T}\mu^N(x,\tau) \Pi_Nu_{\theta}(x,\tau)dxd\tau \\ \notag
&=&\int_{\Omega_T}\left(\nabla u^N\cdot \nabla\Pi_N u_{\theta}(x,\tau)  +q'(u^N) \Pi_N u_{\theta}(x,\tau)+\Pi_N u_{\theta}(x,\tau)(-\Delta)^{\frac{1}{2}}u^N\right) dxd \tau \\ \notag
&=&\int_{\Omega_T}\left(\nabla u^N\cdot \nabla u_{\theta}(x,\tau)  +q'(u^N) \Pi_N u_{\theta}(x,\tau)+u_\theta(-\Delta)^{\frac{1}{2}}u^N\right)dxd \tau.
\end{eqnarray}
Since $\Pi_Nu_{\theta} \rightarrow u_{\theta}$ strongly in $L^2(\Omega_T)$, $\mu^N \rightharpoonup \mu_{\theta}$  in $L^2(\Omega_T)$, by \eqref{eqn:uNwC},\eqref{eqn:quNwC}, passing to the limit in  \eqref{eqn:muNuN} yields
\begin{equation}
\int_{\Omega_T}\mu_{\theta}u_{\theta}dxd\tau=\int_{\Omega_T}\left(|\nabla u_{\theta}|^2+q'(u_{\theta}))u_{\theta}+u_{\theta}(-\Delta)^{\frac{1}{2}}u_{\theta}\right)dxd\tau. \label{eqn:eid}
\end{equation}
\eqref{eqn:guNC} and \eqref{eqn:eid} gives
\begin{equation}
\lim_{N\rightarrow \infty}\int_{\Omega_T}|\nabla u^N|^2 dxd\tau=\int_{\Omega_T}|\nabla u_{\theta}|^2 dxd\tau. \label{eqn:normC}
\end{equation}
 By \eqref{eqn:uNbd}, $\nabla u^N \rightharpoonup \nabla u_{\theta}$ weakly in $L^2(\Omega_T)$, thus \eqref{eqn:normC} implies
\begin{equation}
\nabla u^N \rightarrow \nabla u_{\theta}\text{ strongly in } L^2(\Omega_T).\label{eqn:guNsC}
\end{equation}

By  \eqref{eqn:muNbd} and the lower bound on $M_{\theta}$, we have
$$
\left \Arrowvert \nabla \frac{\mu^N}{g_{\theta}(u^N)}\right \Arrowvert_{L^2(\Omega_T)} \leq C \theta^{-\frac{m}{2}}.
$$

By \eqref{eqn:muN}, \eqref{eqn:uNbd} and \eqref{eqn:qunbd}we have
\begin{eqnarray}\notag
&& \left \arrowvert \int_{\Omega}\frac{\mu^N\phi_1}{g_{\theta}(u_N)}\right \arrowvert dx=\left \arrowvert \int_{\Omega}\mu^N \Pi_N\left(\frac{\phi_1}{g_{\theta}(u^N)}\right)\right \arrowvert dx\\ \label{eqn:muNg}
&\leq &\left \arrowvert \int_{\Omega}\nabla u^N\cdot \nabla \Pi_N\left (\frac{\phi_1}{g_{\theta}(u^N)}\right)dx+\int_{\Omega}q'(u^N)\Pi_N\left(\frac{\phi_1}{g_{\theta}(u^N)}\right)dx\right \arrowvert\\ \notag
&&+\left\arrowvert \int_{\Omega} (-\Delta)^{\frac{1}{2}}u^N \Pi_N\left(\frac{\phi_1}{g_{\theta}(u^N)}\right)dx \right \arrowvert\\ \notag
&=&\left\arrowvert\int_{\Omega}\nabla u^N\cdot \nabla \left (\frac{\phi_1}{g_{\theta}(u^N)}\right)dx+\int_{\Omega}q'(u^N)\Pi_N\left(\frac{\phi_1}{g_{\theta}(u^N)}\right)dx\right \arrowvert \\ \notag
&&+\left\arrowvert\int_{\Omega} (-\Delta)^{\frac{1}{2}}u^N \left(\frac{\phi_1}{g_{\theta}(u^N)}\right)dx \right\arrowvert \\ \notag
&\leq& C\theta^{-m-1}\left \Arrowvert \nabla u^N\right \Arrowvert_{L^2(\Omega)}^2+C\theta^{-m}\left \Arrowvert q'(u^N)\right\Arrowvert_{L^2(\Omega)}\left \Arrowvert \phi_1\right \Arrowvert_{L^2(\Omega)}\\ \notag
&&+C\theta^{-m}\left\Arrowvert \nabla u^N\right\Arrowvert_{L^2(\Omega)}\left\Arrowvert \phi_1\right\Arrowvert_{L^2(\Omega)}\\ \notag
&\leq& C\theta^{-m-1}.
\end{eqnarray}
 Poincare's inequality yields
$$
\left \Arrowvert \frac{\mu^N}{g_{\theta}(u^N)} \right \Arrowvert_{L^2(0,T;H^1(\Omega))} \leq C(\theta^{-m-1}+1).
$$
Thus there exists a $w_{\theta} \in L^2(0,T;H^1(\Omega))$ and a subsequence of $\frac{\mu^N}{g_{\theta}(u^N)}$, not relabeled, such that
\begin{equation}
\frac{\mu^N}{g_{\theta}(u^N)}\rightharpoonup w_{\theta} \text{ weakly in } L^2(0,T;H^1(\Omega)). \label{eqn:mugNC}
\end{equation}
Therefore by \eqref{eqn:MuNC}, \eqref{eqn:mugNC} and Sobolev embedding theorem, we have
\begin{equation}
\mu^N=g_{\theta}(u^N)\cdot \frac{\mu^N}{g_{\theta}(u^N)}\rightharpoonup \mu_{\theta}=g_{\theta}(u_\theta)w_{\theta} \text{ weakly in } L^2(0,T;W^{1,s}(\Omega))\label{eqn:muNC}
\end{equation}
for any $1\leq s<2$. Combining \eqref{eqn:sqMC}, \eqref{eqn:mugNC}and  \eqref{eqn:muNC}, we have
\begin{equation}
\sqrt{M_{\theta}(u^N)}\nabla \frac{\mu^N}{g_{\theta}(u^N)} \rightharpoonup \sqrt{M_{\theta}(u_\theta)}\nabla \frac{\mu_\theta}{g_{\theta}(u_\theta)} \text{ weakly in } L^2(0,T;L^q(\Omega)) \label{eqn:nonlinearwC}
\end{equation}
for any $1\leq q <2$. By \eqref{eqn:muNbd}, we can improve this convergence to
\begin{equation}
\sqrt{M_{\theta}(u^N)}\nabla \frac{\mu^N}{g_{\theta}(u^N)} \rightharpoonup \sqrt{M_{\theta}(u_\theta)}\nabla \frac{\mu_\theta}{g_{\theta}(u_\theta)} \text{ weakly in } L^2(0,T;L^2(\Omega)). \label{eqn:muNl2wC}
\end{equation}

Since $g_{\theta}\geq \theta^{m}$,  \eqref{eqn:uNp} implies
\begin{equation}
\frac{g'(u^N)}{g_{\theta}^{\frac{3}{2}}(u^N)}\rightarrow \frac{g'_{\theta}(u_{\theta})}{g_{\theta}^{\frac{3}{2}}(u_{\theta})} \text{ a.e  in  } \Omega_T. \label{eqn:lngC}
\end{equation}
In addition,  $\frac{g'(u^N)}{g_{\theta}^{\frac{3}{2}}(u^N)}$ is bounded by
\begin{equation}
\left \arrowvert \frac{g'(u^N)}{g_{\theta}^{\frac{3}{2}}(u^N)}\right \arrowvert \leq C \theta^{-1-\frac{m}{2}}. \label{eqn:lngbd}
\end{equation}

It follows from \eqref{eqn:guNsC}, \eqref{eqn:lngC}, \eqref{eqn:lngbd}  and generalized dominated convergence theorem (see Remark \ref{rmk:gdc}) that
\begin{equation}
\frac{g'(u^N)}{g_{\theta}^{\frac{3}{2}}(u^N)}\nabla u^N\rightarrow \frac{g'_{\theta}(u_{\theta})}{g_{\theta}^{\frac{3}{2}}(u_{\theta})}\nabla u_{\theta} \text{ strongly in } L^2(\Omega_T).  \label{eqn:l2conv}
\end{equation}
Let
$$
f^N(t)=\left\Arrowvert\frac{g'(u^N(x,t))}{g_{\theta}^{\frac{3}{2}}(u^N(x,t))}\nabla u^N(x,t)- \frac{g'_{\theta}(u_{\theta}(x,t))}{g_{\theta}^{\frac{3}{2}}(u_{\theta}(x,t))}\nabla u_{\theta}(x,t) \right \Arrowvert_{L^2(\Omega)},
$$
by \eqref{eqn:l2conv}, we can extract  a subsequence of $f^N$, not relabeled, such that  $f^N(t)\rightarrow 0$ a.e. in (0,T). By Egorov's theorem, for any given  $\delta>0$, there exists $T_{\delta}\subset [0,T]$ with $|T_{\delta}|<\delta$ such that $f^N(t)$ converges to $0$ uniformly on $[0,T]\backslash T_{\delta}$.

Given  $\alpha(t) \in L^2(0,T)$, for any $\varepsilon>0$, there exists $T_{\delta}\subset [0,T]$ with $|T_{\delta}|<\delta$ such that
\begin{equation}
\int_{T_\delta}\alpha^2(t) dt <\epsilon. \label{eqn:alphat}
\end{equation}
Multiplying \eqref{eqn:uN} by $\alpha(t)$ and integrating  in time yield
\begin{eqnarray}  \label{eqn:uNl}
&&\int_0^T\alpha(t)\int_{\Omega}\partial_tu^N\phi_j dxdt\\ \notag
&=&-\beta\int_{\Omega_T}\alpha(t)\mu^N\phi_j dxdt-\int_{\Omega_T}\alpha(t)M_{\theta}(u^N)\nabla \frac{\mu^N}{g_{\theta}(u^N)}\cdot \nabla \frac{\phi_j}{g_{\theta}(u^N)}dxdt \\ \notag
&=&-\beta \int_{\Omega_T}\mu^N\alpha(t)\phi_j dxdt-\int_{\Omega_T}M_0 \alpha(t) \nabla \frac{\mu^N}{g_{\theta}(u^N)}\cdot \nabla \phi_j dxdt \\ \notag
&&+\int_{\Omega_T}\alpha(t) \sqrt{M_0}\phi_j \frac{g'_{\theta}(u^N)}{g^{\frac{3}{2}}_{\theta}(u^N)}\nabla u^N\cdot \sqrt{M_{\theta}(u^N)}\nabla \frac{\mu^N}{g_{\theta}(u^N)} dxdt\\ \notag
&&=-A^N-I^N+II^N.
\end{eqnarray}
Since $\alpha(t) \phi_j \in L^2(0,T;H^1(\Omega))$, by \eqref{eqn:mugNC} and \eqref{eqn:muNC}, we have
\begin{equation}
A^N=\beta\int_{\Omega_T}\mu^N \alpha(t)\phi_j dxdt \rightarrow \beta\int_{\Omega_T}\mu_{\theta}\alpha(t)\phi_j dxdt, \label{eqn:Alimit}
\end{equation}
and
\begin{equation}
I^N=\int_{\Omega_T}M_0 \alpha(t) \nabla \frac{\mu^N}{g_{\theta}(u^N)}\cdot \nabla \phi_j dxdt \rightarrow \int_{\Omega_T}M_0 \alpha(t) \nabla \frac{\mu_{\theta}}{g_{\theta}(u_{\theta})}\cdot \nabla \phi_j dxdt. \label{eqn:Ilimit}
\end{equation}
To find the limit of $II^N$, since

\begin{eqnarray}\label{eqn:IIlimit}
&&\\ \notag
&&\int_{\Omega_T}\alpha(t)\phi_j \left(\frac{g'_{\theta}(u^N)}{g^{\frac{3}{2}}_{\theta}(u^N)}\nabla u^N \sqrt{M_{\theta}(u^N)}\nabla \frac{\mu^N}{g_{\theta}(u^N)}- \frac{g'_{\theta}(u_\theta)}{g^{\frac{3}{2}}_{\theta}(u_\theta)}\nabla u_\theta \sqrt{M_{\theta}(u_{\theta})}\nabla \frac{\mu_{\theta}}{g_{\theta}(u_{\theta})}\right)\\ \notag
&=&\int_{\Omega_T}\alpha(t)\phi_j \left(\frac{g'_{\theta}(u^N)}{g^{\frac{3}{2}}_{\theta}(u^N)}\nabla u^N-\frac{g'_{\theta}(u_\theta)}{g^{\frac{3}{2}}_{\theta}(u_\theta)}\nabla u_\theta\right)\cdot \sqrt{M_{\theta}(u^N)}\nabla \frac{\mu^N}{g_{\theta}(u^N)} dxdt\\ \notag
&&+\int_{\Omega_T}\alpha(t)\phi_j \frac{g'_{\theta}(u_\theta)}{g^{\frac{3}{2}}_{\theta}(u_\theta)}\nabla u_\theta\cdot \left( \sqrt{M_{\theta}(u^N)}\nabla \frac{\mu^N}{g_{\theta}(u^N)}-\sqrt{M_{\theta}(u_\theta)}\nabla \frac{\mu_\theta}{g_{\theta}(u_\theta)} \right)dxdt\\ \notag
&=&II^N_1+II^N_2
\end{eqnarray}
From bound
\begin{eqnarray*}
&&\int_{\Omega_T}\left\arrowvert \alpha(t)\phi_j \frac{g'_{\theta}(u_\theta)}{g^{\frac{3}{2}}_{\theta}(u_\theta)}\nabla u_\theta\right\arrowvert^2 dxdt \\
&\leq& C\theta^{-2-m} \left\Arrowvert \nabla u_{\theta}\right \Arrowvert_{L^{\infty}(0,T;L^2(\Omega))}^2 \int_0^T\alpha^2(t)^2 dt,
\end{eqnarray*}
 we conclude that $\alpha(t)\phi_j \frac{g'_{\theta}(u_\theta)}{g^{\frac{3}{2}}_{\theta}(u_\theta)}\nabla u_\theta \in L^2(\Omega_T)$. By \eqref{eqn:muNl2wC}, we can pass to the limit in $II^N_2$ and conclude
$$
II^N_2=\int_{\Omega_T}\alpha(t)\phi_j \frac{g'_{\theta}(u_\theta)}{g^{\frac{3}{2}}_{\theta}(u_\theta)}\nabla u_\theta\cdot \left( \sqrt{M_{\theta}(u^N)}\nabla \frac{\mu^N}{g_{\theta}(u^N)}-\sqrt{M_{\theta}(u_\theta)}\nabla \frac{\mu_\theta}{g_{\theta}(u_\theta)} \right)dxdt\rightarrow 0.
$$
To pass to the limit in $II^N_1$, we write
\begin{eqnarray*}
II^N_1&=&\int_{\Omega_T}\alpha(t)\phi_j \left(\frac{g'_{\theta}(u^N)}{g^{\frac{3}{2}}_{\theta}(u^N)}\nabla u^N-\frac{g'_{\theta}(u_\theta)}{g^{\frac{3}{2}}_{\theta}(u_\theta)}\nabla u_\theta\right)\cdot \sqrt{M_{\theta}(u^N)}\nabla \frac{\mu^N}{g_{\theta}(u^N)} dxdt\\
&=& \int_{T_{\delta}}\int_{\Omega}\alpha(t)\phi_j \left(\frac{g'_{\theta}(u^N)}{g^{\frac{3}{2}}_{\theta}(u^N)}\nabla u^N-\frac{g'_{\theta}(u_\theta)}{g^{\frac{3}{2}}_{\theta}(u_\theta)}\nabla u_\theta\right)\cdot \sqrt{M_{\theta}(u^N)}\nabla \frac{\mu^N}{g_{\theta}(u^N)} dxdt\\
&&+\int_{[0,T]\backslash T_{\delta}}\int_{\Omega}\alpha(t)\phi_j \left(\frac{g'_{\theta}(u^N)}{g^{\frac{3}{2}}_{\theta}(u^N)}\nabla u^N-\frac{g'_{\theta}(u_\theta)}{g^{\frac{3}{2}}_{\theta}(u_\theta)}\nabla u_\theta\right)\cdot \sqrt{M_{\theta}(u^N)}\nabla \frac{\mu^N}{g_{\theta}(u^N)} dxdt\\
&=&II^N_{11}+II^N_{12}.
\end{eqnarray*}

We bound $II^N_{11}$ by
\begin{eqnarray*}
|II^N_{11}|&\leq &\int_{T_{\delta}}|\alpha(t)| \left \Arrowvert\frac{g'_{\theta}(u^N)}{g^{\frac{3}{2}}_{\theta}(u^N)}\nabla u^N-\frac{g'_{\theta}(u_\theta)}{g^{\frac{3}{2}}_{\theta}(u_\theta)}\nabla u_\theta \right \Arrowvert_{L^2(\Omega)}\left\Arrowvert\sqrt{M_{\theta}(u^N)}\nabla \frac{\mu^N}{g_{\theta}(u^N)} \right\Arrowvert_{L^2(\Omega)} dt \\
&\leq& \left \Arrowvert \alpha(t)\right \Arrowvert_{L^2(T_{\delta})} \left \Arrowvert\sqrt{M_{\theta}(u^N)}\nabla \frac{\mu^N}{g_{\theta}(u^N)} \right \Arrowvert_{L^2(\Omega_T)}\left\Arrowvert \frac{g'_{\theta}(u^N)}{g^{\frac{3}{2}}_{\theta}(u^N)}\nabla u^N-\frac{g'_{\theta}(u_\theta)}{g^{\frac{3}{2}}_{\theta}(u_\theta)}\nabla u_\theta \right\Arrowvert_{L^{\infty}(0,T;L^2(\Omega))}\\
&\leq &C(\theta) \varepsilon.
\end{eqnarray*}
For $II^N_{12}$, we have
\begin{eqnarray*}
|II^N_{12}|&\leq &\int_{[0,T]\backslash T_{\delta}}|\alpha(t)| \left \Arrowvert\frac{g'_{\theta}(u^N)}{g^{\frac{3}{2}}_{\theta}(u^N)}\nabla u^N-\frac{g'_{\theta}(u_\theta)}{g^{\frac{3}{2}}_{\theta}(u_\theta)}\nabla u_\theta \right \Arrowvert_{L^2(\Omega)}\left\Arrowvert\sqrt{M_{\theta}(u^N)}\nabla \frac{\mu^N}{g_{\theta}(u^N)} \right\Arrowvert_{L^2(\Omega)} dt\\
&=&\int_{[0,T]\backslash T_{\delta}}|\alpha(t)|f^N(t)|\left\Arrowvert\sqrt{M_{\theta}(u^N)}\nabla \frac{\mu^N}{g_{\theta}(u^N)} \right\Arrowvert_{L^2(\Omega)} dt.
\end{eqnarray*}
Since  $f^N(t)$ converges to $0$ uniformly, $\alpha(t) \in L^2(0,T)$ and $\left\Arrowvert\sqrt{M_{\theta}(u^N)}\nabla \frac{\mu^N}{g_{\theta}(u^N)} \right\Arrowvert_{L^2(\Omega_T)}\leq C $, letting $N \rightarrow \infty $ in $II^N_{12}$ yields $II^N_{12} \rightarrow 0$. Letting $\varepsilon \rightarrow 0$, we conclude $II^N_1\rightarrow 0$ as $N\rightarrow \infty$.
Passing to the limit in \eqref{eqn:uNl}, we have
\begin{eqnarray}\notag
&&\int_0^T\alpha (t) \int_{\Omega}\left <\partial_tu_{\theta},\phi_j \right>_{(W^{1,q}(\Omega))', W^{1,q}(\Omega))} dt\\ \label{eqn:phij}
&=&-\beta\int_{\Omega_T}\alpha(t)\mu_{\theta}\phi_j dxdt-\int_{\Omega_T}\alpha(t)M_{\theta}(u_{\theta})\nabla \frac{\mu_\theta}{g_{\theta}(u_\theta)}\cdot \nabla \frac{\phi_j}{g_\theta(u_\theta)} dxdt.
\end{eqnarray}
Fix $q>2$,  given any $\phi \in L^2(0,T;W^{1,q}(\Omega))$, its Fourier series $\sum_{j=1}^\infty a_j(t) \phi_j(x)$ converges strongly to $\phi$ in $L^2(0,T;W^{1,q}(\Omega))$. Hence
\begin{eqnarray*}
&&\int_{\Omega_T}M_{\theta}(u_{\theta})\nabla \frac{\mu_{\theta}}{g_{\theta}(u_{\theta})}\cdot \nabla \frac{\phi-\Pi_N \phi}{g_{\theta}(u_{\theta})} dxdt \\
&=&\int_{\Omega_T}M_0\nabla \frac{\mu_{\theta}}{g_{\theta}(u_{\theta})}\cdot \nabla (\phi-\Pi_N \phi)dxdt \\
&&-\int_{\Omega_T}(\phi-\Pi_N \phi)\sqrt{M_0}\frac{g_{\theta}'(u_{\theta})}{g_{\theta}^{ \frac{3}{2}}(u_{\theta})}\nabla u_{\theta}\cdot \sqrt{M_{\theta}(u_{\theta})}\nabla \frac{\mu_{\theta}}{g_{\theta}(u_{\theta})} dxdt\\
&=&J^N_1-J^N_2,
\end{eqnarray*}
where  by\eqref{eqn:mugNC}, \eqref{eqn:muNC} and  strong convergence of $\Pi_N\phi $ to $\phi$ in $L^2(0,T;H^1(\Omega))$, we conclude  $$
J^N_1=\int_{\Omega_T}M_0\nabla \frac{\mu_{\theta}}{g_{\theta}(u_{\theta})}\cdot \nabla (\phi-\Pi_N \phi)dxdt \rightarrow 0
$$We can bound $J^N_2$ by
\begin{eqnarray*}
&&|J^N_2|=\left \arrowvert \int_{\Omega_T}(\phi-\Pi_N \phi)\sqrt{M_0}\frac{g_{\theta}'(u_{\theta})}{g_{\theta}^{3/2}(u_{\theta})}\nabla u_{\theta}\cdot \sqrt{M_{\theta}(u_{\theta})}\nabla \frac{\mu_{\theta}}{g_{\theta}(u_{\theta})} dxdt\right \arrowvert  \\
&\leq &\sqrt{M_0}\int_0^T\left\Arrowvert\phi-\Pi_N \phi\right\Arrowvert_{L^{\infty}(\Omega)} \left\Arrowvert\frac{g_{\theta}'(u_{\theta})}{g_{\theta}^{3/2}(u_{\theta})}\nabla u_{\theta}\right\Arrowvert_{L^2(\Omega)} \left\Arrowvert\sqrt{M_{\theta}(u_{\theta})}\nabla \frac{\mu_{\theta}}{g_{\theta}(u_{\theta})}\right\Arrowvert_{L^2(\Omega)}\\
&\leq&\sqrt{M_0}\left\Arrowvert\frac{g_{\theta}'(u_{\theta})}{g_{\theta}^{3/2}(u_{\theta})}\nabla u_{\theta}\right\Arrowvert_{L^{\infty}(0,T;L^2(\Omega))}\left\Arrowvert\sqrt{M_{\theta}(u_{\theta})}\nabla \frac{\mu_{\theta}}{g_{\theta}(u_{\theta})}\right\Arrowvert_{L^2(\Omega_T)}\left\Arrowvert\phi-\Pi_N \phi\right\Arrowvert_{L^2(0,T;W^{1,q}(\Omega))}\\
&\rightarrow& 0 \text{ as } N\rightarrow \infty.
\end{eqnarray*}
Consequently \eqref{eqn:phij} implies
\begin{eqnarray}\notag
&&\int_0^T\left <\partial_tu_{\theta},\phi \right>_{(W^{1,q}(\Omega))', W^{1,q}(\Omega))} dt\\   \label{eqn:uthetaeq}
&=&-\beta \int_{\Omega_T}\mu_{\theta}\phi dxdt -\int_{\Omega_T}M_{\theta}(u_{\theta})\nabla \frac{\mu_\theta}{g_{\theta}(u_\theta)}\cdot \nabla \frac{\phi}{g_\theta(u_\theta)} dxdt
\end{eqnarray}
for all $\phi \in L^2(0,T;W^{1,q}(\Omega))$ with $q>2$. Moreover, since $u^N(x,0)=\Pi_Nu_0(x) \rightarrow u_0(x)$ in $H^1(\Omega)$, we see that $u_\theta(x,0)=u_0(x)$ by \eqref{eqn:uNC}.

\begin{remark}\label{rmk:gdc}
(Generalized dominated convergence theorem) Assume  $E\subset \R^n$ is  measurable. $g_n\rightarrow g $ strongly in $L^q(E)$ for $1\leq q<\infty$ and $f_n$, $f$: $E\rightarrow \R^n$ are measurable functions satisfying
$$
f_n\rightarrow f \text{ a.e.  in } E; \hspace{0.1in} |f_n|^p \leq |g_n|^q \text{ a.e. in } E
$$
with $1\leq p< \infty$, then $f_n\rightarrow f$ in $L^p(E)$.
\end{remark}

\subsubsection{Regularity of $u_{\theta}$} \label{sec:regularity-u theta}
We now consider the regularity of $u_\theta$. Given any $a_j(t) \in L^2(0,T)$, $a_j(t)\phi_j \in L^2(0,T;C(\overline\Omega)$). Integrating   \eqref{eqn:muN} from $0$ to $T$, by \eqref{eqn:quNwC},\eqref{eqn:muNC} and  \eqref{eqn:guNsC}, we have
\begin{eqnarray*}
&&\int_{\Omega_T}\mu_\theta(x,t) a_j(t)\phi_j(x)dxdt\\
&=&\int_{\Omega_T}\left(\nabla u_\theta\cdot a_j(t)\nabla \phi_j+q'(u_\theta)a_j(t)\phi_j+a_j(t)\phi_j(-\Delta)^{\frac{1}{2}}u_\theta\right) dxdt
\end{eqnarray*}
for all $j \in \N$. Given any $\phi \in L^2(0,T;H^1(\Omega))$, its Fouirier series strongly converges to $\phi$ in $L^2(0,T;H^1(\Omega))$, therefore
\begin{eqnarray}\label{eqn:mutheta}
&&\int_{\Omega_T}\mu_\theta(x,t) \phi(x)dxdt
=\int_{\Omega_T}\left(\nabla u_\theta\cdot \nabla \phi+q'(u_\theta)\phi+\phi (-\Delta)^{\frac{1}{2}}u_\theta\right)dxdt.
\end{eqnarray}
Recall $\mu_{\theta} \in L^2(0,T;L^p(\Omega))$ and $q'(u_{\theta})\in L^{\infty}(0,T;L^p(\Omega))$ for any $1\leq p<\infty$, regularity theory implies  $u_\theta \in L^2(0,T;H^2(\Omega))$. Hence
\begin{equation}
\mu_\theta=-\Delta u_\theta+q'(u_\theta)+(-\Delta)^{\frac{1}{2}}u_\theta\text{ a.e. in } \Omega_T.
\end{equation}
By Sobolev embedding theorem, $u_{\theta}\in L^{\infty}(0,T;H^1(\Omega)) \hookrightarrow  L^{\infty}(0,T;L^p(\Omega))$ for any $1\leq p < \infty$. Since growth assumption on $q$ implies $|q''(u)|\leq C(1+|u|^{r-1})$,  pick $p>2$, we have
\begin{eqnarray*}
&&\int_{\Omega}|\nabla q'(u_\theta)|^2 dx=\int_{\Omega}|q''(u_\theta)|^2|\nabla u_\theta|^2 dx \\
&\leq&\left \Arrowvert q''(u_\theta)\right \Arrowvert_{L^{\frac{2p}{p-2}}(\Omega)}^2\left \Arrowvert \nabla u_{\theta}\right \Arrowvert_{L^p(\Omega)}^2\\
&\leq &C\left(1+\left \Arrowvert u_\theta \right \Arrowvert_{L^{\frac{2p}{p-2}(r-1)}(\Omega)}^{2(r-1)}\right)\left \Arrowvert \nabla u_{\theta}\right \Arrowvert_{L^p(\Omega)}^2 \\
&\leq & C\left (1+\left \Arrowvert u_{\theta}\right \Arrowvert^{2(r-1)}_{L^{\infty}(0,T;H^1(\Omega))}\right )\left \Arrowvert \nabla u_{\theta}\right \Arrowvert_{L^p(\Omega)}^2\\
&\leq & C\left (1+\left \Arrowvert u_{\theta}\right \Arrowvert^{2(r-1)}_{L^{\infty}(0,T;H^1(\Omega))}\right )\left \Arrowvert \nabla u_{\theta}\right \Arrowvert_{H^1(\Omega)}^2.
\end{eqnarray*}
Therefore $\nabla q'(u_\theta)=q''(u_\theta)\nabla u_\theta \in L^2(\Omega_T)$ with
\begin{eqnarray*}
\int_{\Omega_T}|\nabla q'(u_\theta)|^2 dxdt\leq C \left (1+\left \Arrowvert u_{\theta}\right \Arrowvert^{2(r-1)}_{L^{\infty}(0,T;H^1(\Omega))}\right )\left \Arrowvert \nabla u_{\theta}\right \Arrowvert_{L^2(0,T;H^1(\Omega))}^2.
\end{eqnarray*}
Hence  $q'(u_\theta) \in L^2(0,T;H^1(\Omega))$, combined with $\mu_{\theta} \in L^2(0,T;W^{1,s}(\Omega))$ for any $1\leq s<2$ and $(-\Delta)^{\frac{1}{2}}u_{\theta} \in L^2(0,T;H^1(\Omega))$, we have $u_\theta \in L^2(0,T;W^{3,s}(\Omega))$ and
\begin{equation}
\nabla \mu_\theta=-\nabla \Delta u_{\theta}+q''(u_\theta)\nabla u_\theta +\nabla (-\Delta)^{\frac{1}{2}}u_\theta \text{ a.e. in }  \Omega_T.
\end{equation}
Regularity of $u_{\theta}$ implies  $\nabla u_{\theta} \in L^{\infty}(0,T;L^2(\Omega))\cap L^2(0,T; L^{\infty}(\Omega))$. A simple interpolation shows  $\nabla u_{\theta} \in L^{\frac{2s}{s-2}}(0,T;L^s(\Omega))$ for any $s>2$.  Given any $\phi \in L^p(0,T;W^{1,q}(\Omega))$ with $p>2$ and $q>2$, we have  $g_{\theta}(u_{\theta})\phi \in L^2(0,T;W^{1.r}(\Omega))$ for any $r<q$. Picking  $g_{\theta}(u_{\theta})\phi$ as a test function in \eqref{eqn:uthetaeq}, we have
\begin{equation}
\int_{\Omega_T}\partial_tu_{\theta} g_{\theta}(u_\theta) \phi dxdt=-\beta\int_{\Omega_T}g_{\theta}(u_{\theta})\mu_{\theta}\phi dxdt-\int_{\Omega_T}M_{\theta}(u_\theta)\nabla \frac{\mu_\theta}{g_{\theta}(u_{\theta})}\cdot \nabla \phi dxdt \label{eqn:utheta2}
\end{equation}
for any $\phi \in L^p(0,T;W^{1,q}(\Omega))$ with $p,q>2$.
\begin{remark}
In fact, since $M_{\theta}(u_{\theta}) \in L^{\infty}(0,T;L^p(\Omega))$ for $1\leq p<\infty$, the right hand side of \eqref{eqn:utheta2} is well defined for any $\phi \in L^2(0,T,W^{1,q}(\Omega))$ and we can extend  \eqref{eqn:utheta2} to hold for all $\phi \in L^2(0,T,W^{1,q}(\Omega))$.
\end{remark}
\subsubsection{Energy Inequality}

Since $u^N$ and $\mu^N$ satisfies energy identity \eqref{eqn:energyid}, passing to the limit as $N\rightarrow \infty$ and using the weak convergence of $u^N$,
$q'(u^N)$ and $\sqrt{M_\theta(u^N) }\nabla \frac{\mu^N}{g_\theta(u^N)}$, the energy inequality \eqref{eqn:engineq} follows.

This finishes the proof of Proposition \ref{thm:ndg}.

\subsection{Phase field model with degenerate mobility}
\label{sec:dg}

In this subsection, we prove Theorem \ref{thm:dg}.

Fix initial data $u_0 \in H^1(\Omega)$. We pick a montone decreasing  positive sequence $\theta_i$ with $\lim_{i\rightarrow \infty}\theta_i=0$. By Proposition \ref{thm:ndg} and \eqref{eqn:utheta2}, for each $\theta_i$, there exists
$$
u_i \in L^{\infty}(0,T;H^1(\Omega))\cap L^2(0,T;W^{3,s}(\Omega))\cap C([0,T];L^p(\Omega))
$$
 with weak derivative
$$
\partial_t u_i \in L^2(0,T; (W^{1,q}(\Omega))'),
$$
where $1\leq p<\infty$, $1\leq s <2$, $q>2$ such that $u_{\theta_i}(x,0)=u_0(x)$ and for all $\phi \in L^2(0,T;W^{1,q}(\Omega))$,
\begin{eqnarray}
\label{eqn:ui}\int_{\Omega_T}\partial_t u_i\phi dxdt&=&-\beta\int_{\Omega_T}\mu_i\phi dxdt -\int_{\Omega_T}M_i(u_i)\nabla \frac{\mu_i}{g_i(u_i)}\nabla \frac{\phi}{g_i(u_i)} dxdt,\\\mu_i&=&-\Delta u_i+q'(u_i)+(-\Delta)^{\frac{1}{2}}u_i.
\end{eqnarray}
Moreover,  for all $\phi \in L^p(0,T;W^{1,q}(\Omega))$ with $p,q >2$, the following holds:
\begin{equation} \label{eqn:ui2}
\int_{\Omega_T}g_i(u_i)\partial_t u_i\phi dxdt=-\beta \int_{\Omega_T}g_i(u_i)\mu_i\phi dxdt-\int_{\Omega_T}M_i(u_i)\nabla \frac{\mu_i}{g_i(u_i)}\nabla \phi dxdt.
\end{equation}
Here we write $u_i=u_{\theta_i}$, $M_i(u_i)=M_{\theta_i}(u_{\theta_i})$, $g_i(u_i)=g_{\theta_i}(u_{\theta_i})$ for simplicity of notations. Noticing the bound in \eqref{eqn:uNbd} and \eqref{eqn:muNbd} only depends on $u_0$, we can find a constant $C$, independent of $\theta_i$ such that
\begin{eqnarray}\label{eqn:uibd}
&&\left\Arrowvert \mu_i\right\Arrowvert_{L^2(\Omega_T)}, \left \Arrowvert u_i\right \Arrowvert_{L^{\infty}(0,T;H^1(\Omega))}\leq C, \\ \label{eqn:Mui}
&& \left \Arrowvert \sqrt{M_i(u_i)}\nabla \frac{\mu_i}{g_i(u_i)}\right \Arrowvert_{L^2(\Omega_T)} \leq C.
\end{eqnarray}
Growth condition on $q'$,  and Sobolev embedding theorem gives
\begin{eqnarray*}
&&\left \Arrowvert q'(u_i)\right \Arrowvert_{L^\infty(0,T;L^p(\Omega))} \leq C,\\
&&\left \Arrowvert M_i(u_i)\right \Arrowvert_{L^{\infty}(0,T;L^p(\Omega))} \leq C
\end{eqnarray*}
for any $1\leq p <\infty$. By \eqref{eqn:ui2}, for any $\phi \in L^p(0,T;W^{1,q}(\Omega))$ with $p,q>2$,
\begin{eqnarray*}
&&\left \arrowvert \Int_{\bs \Omega_T}g_i(u_i)\partial_t u_i\phi dxdt\right \arrowvert=\left \arrowvert \Int_{\bs\Omega_T}\left[\beta g_i(u_i)\mu_i \phi +M_i(u_i)\nabla \frac{\mu_i}{g_i(u_i)}\nabla \phi \right]dxdt\right \arrowvert \\
&\leq& \beta \Int_{\bs 0}^T \left(\left \Arrowvert \mu_i\right \Arrowvert_{L^2(\Omega)}\left \Arrowvert g_i(u_i)\right \Arrowvert_{L^{\frac{2q}{q-2}}(\Omega)}\left \Arrowvert \phi\right \Arrowvert_{L^q(\Omega)}\right)dt \\
&&+\Int_{\bs 0}^T\left(\left \Arrowvert \sqrt{M_i(u_i)}\nabla \frac{\mu_i}{g_i(u_i)}\right \Arrowvert_{L^2(\Omega)}\left \Arrowvert \sqrt{M_i(u_i)}\right \Arrowvert_{L^{\frac{2q}{q-2}}(\Omega)}\left \Arrowvert \nabla \phi \right \Arrowvert_{L^q(\Omega)} \right )dt\\
&\leq&\beta\left \Arrowvert g_i(u_i)\right \Arrowvert_{L^{\frac{2p}{p-2}}(0,T;L^{\frac{2q}{q-2}}(\Omega))}\left \Arrowvert \mu_i \right \Arrowvert_{L^2(\Omega_T)}\left \Arrowvert \phi \right \Arrowvert_{L^p(0,T;L^q(\Omega))} \\
&&+\left \Arrowvert M_i(u_i)\right \Arrowvert^{\frac{1}{2}}_{L^{\frac{p}{p-2}}(0,T;L^{\frac{q}{q-2}}(\Omega))} \left \Arrowvert \sqrt{M_i(u_i)}\nabla \frac{\mu_i}{g_i(u_i)}\right \Arrowvert_{L^2(\Omega_T)}\left \Arrowvert \nabla \phi \right \Arrowvert_{L^p(0,T;L^q(\Omega))} \\
&\leq & C  \left \Arrowvert \phi \right \Arrowvert_{L^p(0,T;W^{1,q}(\Omega))}.
\end{eqnarray*}
Let
\begin{equation}
G_i(u_i)=\int_0^{u_i} g_i(a) da. \label{eqn:Gidef}
\end{equation}  Then  $\partial_t G_i(u_i)=g_i(u_i)\partial_t u_i \in L^{p'}(0,T; (W^{1,q}(\Omega))')$ with $p'=\frac{p}{p-1}$ and
\begin{equation}
\left \Arrowvert \partial_t G_i(u_i) \right \Arrowvert_{L^{p'}(0,T; (W^{1,q}(\Omega))')} \leq C \text { for all } i. \label{eqn:Git}
\end{equation}
Moreover, by growth assumption on $g$ and estimates on $u_i$, we have
 \begin{equation}
\left \Arrowvert G_i(u_i)\right\Arrowvert_{ L^{\infty}(0,T;W^{1,s}(\Omega))} \leq C. \label{eqn:Gibd}
\end{equation}
 for $1\leq s <2$.
By  \eqref{eqn:uibd}, \eqref{eqn:Mui}-\eqref{eqn:Gibd} and Remark \ref{rmk:cpt} we can  find a subsequence, not relabeled, a function $u\in L^{\infty}(0,T;H^1(\Omega))$, a function $\mu \in L^2(\Omega_T)$,  a function $\xi \in L^2(\Omega_T)$  and a function $\eta\in L^{\infty}(0,T;W^{1,s}(\Omega)) $ such that as $i \rightarrow \infty$,
\begin{eqnarray}\label{eqn:uiwC}
&&u_i \rightharpoonup u \text{ weakly-* in } L^{\infty}(0,T;H^1(\Omega)), \\ \label{eqn:muiwC}
&& \mu_i\rightharpoonup \mu \text{ weakly in } L^2(\Omega_T), \\
\label{eqn:uip}
&&\sqrt{M_i(u_i)}\nabla \frac{\mu_i}{g_i(u_i)} \rightharpoonup \xi \text{ weakly in } L^2(\Omega_T),\\
&&G_i(u_i)\rightharpoonup \eta\text{ weakly-* in } L^{\infty}(0,T;W^{1,s}(\Omega))\\\label{eqn:GiuiSCp}
&&G_i(u_i) \rightarrow \eta \text{ strongly in } L^{\alpha}(0,T; L^{\beta}(\Omega)) \text{ and a.e. in } \Omega_T,\\ \label{eqn:GuiSCinf}
&&G_i(u_i) \rightarrow \eta \text{ strongly in } C(0,T; L^{\beta}(\Omega)),\\ \label{eqn:devGi}
&& \partial_t G_i(u_i) \rightharpoonup \partial_t \eta \text{ weakly in } L^{p'}(0,T;(W^{1,q}\Omega))').
\end{eqnarray}
where  $1\leq \alpha,\beta <\infty$. By \eqref{eqn:GuiSCinf} and \eqref{eqn:unicov} from Remark \ref{rmk:cptns}, we have
$$
\left \Arrowvert G_i(u_i(x,t+h))-G_i(u_i(x,t))\right\Arrowvert_{C([0,T];L^{\beta}(\Omega))} \rightarrow 0 \text{ uniformly in } i \text{ as } h\rightarrow 0.
$$ Thus given any $\e >0$, there exists $h_{\e}>0$ such that for all $0<h<h_{\e}$ and all $i$,

$$
\left \Arrowvert G_i(u_i(x,t+h))-G_i(u_i(x,t))\right\Arrowvert_{C([0,T];L^{\beta}(\Omega))}^{\beta} <\e.
$$
Given any $\delta>0$, let $I_{\delta}=(1-\delta,1+\delta)\cup (-1-\delta,-1+\delta)$. Consider  the interval having  $u_i(x,t)$ and $u_i(x,t+h)$ as end points.  Denote this interval by $J_i(x,t;h)$. We consider three cases.

\textbf{ Case I: $J_i(x,t;h)\cap I_{\delta}=\varnothing$.}

In this case, $g_i(s)\geq \max\{\theta_i^m, \delta^m\}$ for any $s\in J_i(x,t;h)$ and
$$
\left |G_i(u_i(x,t+h))-G_i(u_i(x,t))\right|=\left|\int_{u_i(x,t)}^{u_i(x,t+h)}g_i(s)ds \right|\geq \delta^m|u_i(x,t+h)-u_i(x,t)|.
$$

\textbf{ Case II: $J_i(x,t;h)\cap I_{\delta}\neq \varnothing$ and $|u_i(x,t+h)-u_i(x,t)|\geq 3\delta$ .}

In this case, we have $$|J_i(x,t;h)\cap I_{\delta}^{c}|\geq \frac{1}{3}|J_i(x,t;h)|$$ and
\begin{eqnarray*}
\left |G_i(u_i(x,t+h))-G_i(u_i(x,t))\right|&\geq& \left|\int_{J_i(x,t;h)\cap I_{\delta}^{c}}g_i(s)ds \right|\\
&\geq& \frac{\delta^m}{3}|u_i(x,t+h)-u_i(x,t)|.
\end{eqnarray*}

\textbf{ Case III: $J_i(x,t;h)\cap I_{\delta}\neq \varnothing$ and $|u_i(x,t+h)-u_i(x,t)|< 3\delta$ }

In this case, we have
$$
g_i(s)\leq \max\{(8\delta+16\delta^2)^m,\theta_i^m\} \text{ for any } s \in J_i(x,t;h).
$$Thus
\begin{eqnarray*}
\left |G_i(u_i(x,t+h))-G_i(u_i(x,t))\right|\leq 3\delta\max\{(8\delta+16\delta^2)^m,\theta_i^m\}.
\end{eqnarray*}

Pick $\delta=\e^{\frac{1}{2m\beta}}$ and fix $t$.  Let
$$
\Omega_i=\{x\in \Omega: J_i(x,t:h) \text{ satisfies case I or II} \}.
$$ Then
\begin{eqnarray*}
&&\int_{\Omega}\left \arrowvert u_i(x,t+h)-u_i(x,t) \right\arrowvert^{\beta} dx \\
&=&\int_{\Omega_{i}}\left \arrowvert u_i(x,t+h)-u_i(x,t) \right\arrowvert^{\beta} dx+\int_{\Omega\backslash\Omega_{i}}\left \arrowvert u_i(x,t+h)-u_i(x,t) \right\arrowvert^{\beta} dx\\
&\leq& 3^{\beta}\e^{-\frac{1}{2}}\int_{\Omega_{i}}\left \arrowvert G_i(u_i(x,t+h))-G_i(u_i(x,t)) \right\arrowvert^{\beta} dx+\int_{\Omega\backslash\Omega_{i}}\left \arrowvert u_i(x,t+h)-u_i(x,t) \right\arrowvert^{\beta} dx\\
&\leq&{  3^{\beta}}\e^{\frac{1}{2}}+C\e^{\frac{1}{2m}}
\end{eqnarray*}
Taking maximum over $t\in [0,T]$ on the left side, we have for all $i$, any $h<h_{\e}$,
$$
\left \Arrowvert (u_i(x,t+h)-u_i(x,t)\right\Arrowvert_{C([0,T];L^{\beta}(\Omega))}^{\beta} < \e^{\frac{1}{2}}+C\e^{\frac{1}{2m}}.
$$

Thus
$$
\left \Arrowvert u_i(x,t+h)-u_i(x,t)\right\Arrowvert_{C([0,T];L^{\beta}(\Omega))}^{\beta} \rightarrow 0 \text{ uniformly as } h\rightarrow 0.
$$

In addition, for any $0<t_1<t_2<T$, \eqref{eqn:uibd} implies
$$
\int_{t_1}^{t_2} u_i(x,t)dt \text{ is relatively compact in }L^{\beta}(\Omega).
$$
Therefore we conclude from Remark \ref{rmk:cptns} that
\begin{equation}
u_i\rightarrow u(x,t) \text{ strongly in  }C([0,T];L^{\beta}(\Omega)) \text{ for } 1\leq \beta <\infty.\label{eqn:uistronginf}
\end{equation}
Similarly. we can prove
\begin{equation}
u_i\rightarrow u(x,t) \text{ strongly in  }L^{\alpha}(0,T;L^{\beta}(\Omega)) \text{ for } 1\leq \alpha,\beta <\infty \text{ and a.e. in }\Omega_T.\label{eqn:uistronglp}
\end{equation}
Growth condition on $M(u)$ and \eqref{eqn:uistronginf}, \eqref{eqn:uistronglp} yield
\begin{eqnarray}\label{eqn:miuicv}
&&M_i(u_i)\rightarrow M(u) \text{ strongly in }C([0,T];L^{\beta}(\Omega)) \text{ for } 1\leq \beta <\infty, \\ \label{eqn:miuicvlp}
&&M_i(u_i)\rightarrow M(u) \text{ strongly in }L^{\alpha}(0,T;L^{\beta}(\Omega))  \text{ for } 1\leq \alpha,\beta <\infty,\\ \label{eqn:sqmicv}&&\sqrt{M_i(u_i)}\rightarrow \sqrt{M(u)} \text{ strongly in  }C([0,T];L^{\gamma}(\Omega)) \text{ for } 1\leq \gamma <\infty.
\end{eqnarray}
Hence $G_i(u_i)$ converges to $G(u)$ a.e. in $\Omega_T$ and $\eta=G(u)$.
Passing to the limit in \eqref{eqn:ui2}, by \eqref{eqn:uibd}, \eqref{eqn:uip}, \eqref{eqn:devGi}, \eqref{eqn:miuicv} and  \eqref{eqn:sqmicv}, we have
\begin{eqnarray}\label{eqn:ueq}
&&\int_0^T\left<g(u)\partial_tu, \phi\right>_{((W^{1,q}(\Omega))',W^{1,q}(\Omega))} dt\\ \notag
&=&-\beta\int_{\Omega_T}g(u)\mu \phi dxdt-\int_{\Omega_T} \sqrt{M(u)}\xi\cdot \nabla \phi dxdt
\end{eqnarray}
for any $\phi \in L^p(0,T;W^{1,q}(\Omega))$ with $p,q>2$.
\begin{remark} \label{rmk:cptns}
(Compactness in $L^p(0,T;B)$ Theorem 1 in \cite{Sim86})
Assume $B$ is a Banach space and $F\subset L^p(0,T;B)$. $F$ is relatively compact in  $L^p(0,T;B)$ for $1\leq p<\infty$, or in $C([0,T],B)$ for $p=\infty$ if and only if
\begin{equation}
\left\{\int_{t_1}^{t_2}f(t)dt: f\in F\right \} \text{ is relatively compact in }B, \forall   0<t_1<t_2<T \label{eqn:unibd}
\end{equation}
\begin{equation}
\left\Arrowvert\tau_hf-f\right \Arrowvert_{L^p(0,T;B)}\rightarrow 0 \text{ as } h\rightarrow 0 \text{ uniformly for } f\in F. \label{eqn:unicov}
\end{equation}
Here $\tau_h f(t)=f(t+h)$ for $h>0$ is defined on $[-h,T-h]$.
\end{remark}

\subsubsection{Weak convergence of  $\nabla \frac{\mu_i}{g_i(u_i)}$}
We now look for relation between $\xi$ and $u$. Following the  idea in \cite{DaiDu16}, we  decompose $\Omega_T$ as follows. Let  $\delta_j$ be a positive sequence monotonically decreasing to $0$. By \eqref{eqn:uip} and Egorov's theorem, for every $\delta_j>0$, there exists $B_j \subset \Omega_T$ satisfying $|\Omega_t\backslash B_j|< \delta_j$ such that
\begin{equation}
u_i\rightarrow u \text{ uniformly in } B_j. \label{eqn:unifui}
\end{equation}
We can pick
\begin{equation}
B_1\subset B_2 \subset \cdots \subset B_j\subset B_{j+1} \subset\cdots \subset \Omega_T. \label{eqn:Bj}
\end{equation}
Define
$$
P_j:=\{(x,t)\in \Omega_T:|1-u^2| > \delta_j\}.
$$
Then
\begin{equation}
P_1\subset P_2 \subset \cdots \subset P_j \subset P_{j+1}\subset \cdots \subset\Omega_T. \label{eqn:Pj}
\end{equation}
Let  $B=\cup _{j=1}^{\infty}B_j$ and $P=\cup_{j=1}^{\infty}P_j$.  Then $|\Omega_T\backslash B|=0$ and each $B_j$ can be split into two parts:
\begin{eqnarray*}
&&D_j=B_j\cap P_j, \text{ where } |1-u^2|>\delta_j, \text{ and } u_i\rightarrow u \text{ uniformly},\\
&&\tilde D_j=B_j\backslash P_j, \text{ where } |1-u^2|\leq \delta_j, \text{ and } u_i \rightarrow u \text { uniformly }.
\end{eqnarray*}
\eqref{eqn:Bj} and \eqref{eqn:Pj} imply
\begin{equation}
D_1\subset D_2 \subset \cdots \subset D_j \subset D_{j+1}\subset \cdots \subset D:=B\cap P. \label{eqn:Dj}
\end{equation}
For any $\Psi \in L^p(0,T;L^q(\Omega,\R^n))$ with $ p,q>2$, we have
\begin{eqnarray}\notag
&&\int_{\Omega_T}M_i(u_i)\nabla \frac{\mu_i}{g_i(u_i)}\cdot \Psi dxdt\\ \notag
&=&\int_{\Omega_T\backslash B_j}M_i(u_i)\nabla \frac{\mu_i}{g_i(u_i)}\cdot \Psi dxdt+\int_{D_j}M_i(u_i)\nabla \frac{\mu_i}{g_i(u_i)}\cdot \Psi dxdt\\ \label{eqn:rhsi}
&&+\int_{\tilde D_j}M_i(u_i)\nabla \frac{\mu_i}{g_i(u_i)}\cdot \Psi dxdt
\end{eqnarray}
The left hand side of \eqref{eqn:rhsi} converges to $\int_{\Omega_T}\sqrt{M(u)}\xi\cdot \Psi dxdt$. We analyze the three terms on the right hand side separately.
 To estimate the first term on the right hand side of \eqref{eqn:rhsi}, noticing $|\Omega_T\backslash B_j| \rightarrow 0$ and
\begin{eqnarray*}
\lim_{i\rightarrow \infty} \int_{\Omega_T\backslash B_j}M_i(u_i)\nabla \frac{\mu_i}{g_i(u_i)}\cdot \Psi dxdt=\ \int_{\Omega_T\backslash B_j}\sqrt{M(u)}\xi\cdot \Psi dxdt,
\end{eqnarray*}
we have
\begin{eqnarray*}
\lim_{j\rightarrow \infty}\lim_{i\rightarrow \infty}\int_{\Omega_T\backslash B_j}M_i(u_i)\nabla \frac{\mu_i}{g_i(u_i)}\cdot \Psi dxdt=0.
\end{eqnarray*}
By uniform convergence of $u_i$ to $u$ in $B_j$, we introduce subsequence $u_{j,k}$ such that $u_{j,k}\rightarrow u$ uniformly in $B_j$ and there exists $N_j $ such that for all $k\geq N_j$,
\begin{equation}
|1-u^2_{j,k}|>\frac{\delta_j}{2} \text{ in } D_j,    \hspace{0.2 in} |1-u^2_{j,k}|\leq 2\delta_j \text{ in } \tilde D_j.
\end{equation}
Thus the third term on the right hand side of \eqref{eqn:rhsi} can be estimated by
\begin{eqnarray*}
&&\lim_{j\rightarrow \infty}\lim_{k\rightarrow \infty}\left \arrowvert \int_{\tilde D_j}M_{j,k}(u_{j,k})\nabla \frac{\mu_{j,k}}{g_{j,k}(u_{j,k})} \cdot \Psi dxdt \right \arrowvert \\
&\leq&\lim_{j\rightarrow \infty}\lim_{k\rightarrow \infty}\left \{\left(\sup_{\tilde D_j}\sqrt{M_{j,k}(u_{j,k})}\right)\left \Arrowvert \Psi \right \Arrowvert_{L^2(\tilde D_j)}\left \Arrowvert \sqrt{M_{j,k}(u_{j,k})}\nabla \frac{\mu_{j,k}}{g_{j,k}(u_{j,k})}\right \Arrowvert_{L^2(\tilde D_j)}\right \} \\
&\leq&\left(\sup_{\tilde D_j}\sqrt{M_{j,k}(u_{j,k})}\right)|\Omega|^{\frac{q-2}{2q}}\left \Arrowvert \Psi \right \Arrowvert_{L^2(0,T;L^q(\Omega)}\left \Arrowvert \sqrt{M_{j,k}(u_{j,k})}\nabla \frac{\mu_{j,k}}{g_{j,k}(u_{j,k})}\right \Arrowvert_{L^2(\tilde D_j)}\\
&&\leq C\lim_{j\rightarrow \infty}\lim_{k\rightarrow \infty}\max\left \{(2\delta_j)^{m/2}, \theta_{j,k}^{m/2}\right\}\\
&=&0.
\end{eqnarray*}
For the second  term, we see that
\begin{eqnarray*}
&&\left(\frac{\delta_j}{2}\right)^m\int_{D_j}|\nabla \frac{\mu_{j,k}}{g_{j,k}(u_{j,k})}|^2 dxdt \\
&\leq&\int_{D_j}M_{j,k}(u_{j,k})|\nabla \frac{\mu_{j,k}}{g_{j,k}(u_{j,k})}|^2 dxdt\\
&\leq& \int_{\Omega_T}M_{j,k}(u_{j,k})|\nabla \frac{\mu_{j,k}}{g_{j,k}(u_{j,k})}|^2 dxdt \leq C.
\end{eqnarray*}
Therefore $\nabla \frac{\mu_{j,k}}{g_{j,k}(u_{j,k})}$ is bounded in $L^2(D_j)$ and we can extract a further subsequence, not relabeled, which converges weakly to some $\xi_j \in L^2(D_j)$.
Since $D_j$ is an increasig sequence of sets with $\lim_{j\rightarrow \infty} D_j=D$, we have $\xi_j=\xi_{j-1}$ a.e. in $D_{j-1}$. By setting $\xi_j=0$ outside $D_j$, we can extend $\xi_j $ to a $L^2$  function $\tilde \xi_j$ defined in $D$. Therefore for a.e. $x\in D$, there exists a limit of $\tilde \xi_j(x)$ as $j \rightarrow \infty$. Let $\xi(x)=\lim_{j\rightarrow \infty}\tilde \xi_j(x)$, we see that $\xi(x)=\xi_j(x)$ for a.e $x \in D_j$ and for all $j$.

By a standard diagonal  argument, we can extract a subsequnce such that
\begin{equation}
\nabla \frac{\mu_{k,N_k}}{g_{k,N_k}(u_{k,N_k})}\rightharpoonup \zeta \text{ weakly in } L^2(D_j) \text{ for all } j. \label{eqn:mugwC}
\end{equation}

By strong convergence of $\sqrt{M_i(u_i)}$ to $\sqrt{M(u)}
$ in $C([0,T];L^{\beta}(\Omega))$ for $1\leq \beta<\infty$, we obtain
$$
\chi_{D_j}\sqrt{M_{k,N_k}(u_{k,N_k})}\nabla \frac{\mu_{k,N_k}}{g_{k,N_k}(u_{k,N_k})}\rightharpoonup \chi_{D_j}\sqrt{M(u)}\zeta
$$
weakly in $L^2(0,T;L^q(\Omega))$ for $1\leq q <2$ and all $j$.  Recall $\sqrt{M_i(u_i)}\nabla \frac{\mu_i}{g_i(u_i)} \rightarrow \xi$ weakly in $L^2(\Omega_T)$, we have $\xi=\sqrt{M(u)}\zeta$ in $D_j$ for all $j$. Hence $\xi=\sqrt{M(u)}\zeta$ in $D$ and consequently
$$
\chi_DM_{k,N_k}(u_{k,N_k})\nabla \frac{\mu_{k,N_k}}{g_{k,N_k}(u_{k,N_k})}\rightharpoonup \chi_{D}M(u)\zeta
$$
weakly in $L^2(0,T;L^q(\Omega))$ for $1\leq q <2$.

Replacing $u_i$ by subsequence $u_{k,N_k}$ in \eqref{eqn:rhsi} and letting $k \rightarrow \infty$ then $j\rightarrow \infty$, we have
\begin{eqnarray}\label{eqn:limit}
\int_{\Omega_T}\sqrt{M(u)}\xi \cdot \Psi dxdt&=&\lim_{j\rightarrow \infty}\int_{D_j}M(u)\zeta \cdot \Psi dxdt \\ \notag
&=&\int_D M(u)\zeta\cdot \Psi dxdt.
\end{eqnarray}
It follows from \eqref{eqn:ueq} and \eqref{eqn:limit} that
\begin{eqnarray}
&&\int_0^T\left<g(u)\partial_t u,\phi\right>_{((W^{1,q}(\Omega))',W^{1,q}(\Omega))} dt\\ \notag
&=&-\beta\int_{\Omega_T}g(u)\mu \phi dxdt-\int_{D}M(u)\zeta\cdot \nabla \phi dxdt
\end{eqnarray}
for all $\phi \in L^p(0,T;W^{1,q}(\Omega))$ where $p,q>2$.

\subsubsection{Relation between $\zeta$ and $u$}

The desired relation between $\zeta$ and $u$ is
\begin{eqnarray}\label{eqn:zeta}
\zeta&=&\frac{1}{g(u)}\nabla \mu-\mu \frac{g'(u)}{g^2(u)}\nabla u \\ \label{eqn:muae}
\mu&=& -\Delta u+q'(u)+ (-\Delta)^{\frac{1}{2}}u.
\end{eqnarray}
Given the known regularity $u\in L^{\infty}(0,T;H^1(\Omega))$ and degeneracy of $g(u)$,  the right hand side of \eqref{eqn:zeta} might not be defined as a function.   We can, however,
under  suitable assumptions on integrability of $\nabla \Delta u$, find an explicit expression of $\zeta$ in {  the form} of \eqref{eqn:zeta}-\eqref{eqn:muae}  in suitable subset of $\Omega_T$.

{\it{Claim I: If  for some $j$, the interior of $D_j$, denoted by $(D_j)^{\circ}$, is not empty, then
$$
\nabla \Delta u \in L^1((D_j)^{\circ})
$$
and
$$
\zeta=\frac{-\nabla \Delta u+q''(u)\nabla u+\nabla (-\Delta)^{\frac{1}{2}}u}{g(u)}
-\frac{g'(u)}{g^2(u)}\left(-\Delta u+q'(u)+(-\Delta)^{\frac{1}{2}} u \right)\nabla u
$$
 a.e. in  $(D_j)^{\circ}$.
}}

{\it{Proof of the claim I.}} Since
\begin{equation}
\mu_{k,N_k}=-\Delta u_{k,N_k}+q'(u_{k,N_k})+(-\Delta)^{\frac{1}{2}} u_{k,N_k} \text{ in } \Omega_T, \label{eqn:muik}
\end{equation}
The right hand side of \eqref{eqn:muik}  converges to $-\Delta u+q'(u)+(-\Delta)^{\frac{1}{2}}u$ in distributional sense while the left side converges weakly to $\mu$ in $L^2(\Omega_T)$. Hence
$$
\mu=-\Delta u+q'(u)+(-\Delta)^{\frac{1}{2}}u \text{ in }  L^2(\Omega_T).
$$
Therefore $u \in L^2(0,T;H^2(\Omega))$.
On the other hand, using  $u_{k,N_k}$ and $u$ as test functions  in  \eqref{eqn:mutheta} yield
\begin{eqnarray*}
&&\int_{\Omega_T}\mu_{k,N_k}u_{k,N_k}dxdt=\int_{\Omega_T}\left (\left \arrowvert \nabla u_{k,N_k}\right\arrowvert^2+q'(u_{k,N_k})u_{k,N_k}+u_{k,N_k}(-\Delta)^{\frac{1}{2}}u_{k,N_k} \right)dxdt\\
&&\int_{\Omega_T}\mu_{k,N_k}udxdt=\int_{\Omega_T}\left (\nabla u_{k,N_k}\cdot \nabla u+q'(u_{k,N_k})u+u(-\Delta)^{\frac{1}{2}}u_{k,N_k} \right)dxdt.
\end{eqnarray*}
Passing to the limit, by \eqref{eqn:uistronglp}, growth assumptions on $q'$ and \eqref{eqn:muiwC}, we have
$$
\lim_{k\rightarrow \infty} \int_{\Omega_T}\left \arrowvert \nabla u_{k,N_k}\right\arrowvert^2=\int_{\Omega_T}\left \arrowvert \nabla u\right\arrowvert^2.
$$
Therefore $$\nabla u_{k,N_k} \rightarrow \nabla u \text{ strongly in } L^2(\Omega_T).$$
Since $u_{k,N_k} \in L^2(0,T;W^{3,s}(\Omega))$,  we can differentiate \eqref{eqn:muik} and get
\begin{equation}
\nabla \mu_{k,N_k}=-\nabla \Delta u_{k,N_k}+q''(u_{k,N_k})\nabla u_{k,N_k}+\nabla (-\Delta)^{\frac{1}{2}} u_{k,N_k}, \label{eqn:gmuk}
\end{equation}
and
\begin{equation}
\nabla \frac{\mu_{k,N_k}}{g_{k,N_k}(u_{k,N_k})}=\frac{1}{g_{k,N_k}(u_{k,N_k})}\nabla \mu_{k,N_k}-\mu_{k,N_k}\frac{g'_{k,N_k}(u_{k,N_k})}{g^2_{k,N_k}(u_{k,N_k})}\nabla u_{k,N_k} \label{eqn:muikuik}
\end{equation}
on $D_j^{\circ}$.
Thus
\begin{equation}
\nabla \mu_{k,N_k}=g_{k,N_k}(u_{k,N_k})\nabla  \frac{\mu_{k,N_k}}{g_{k,N_k}(u_{k,N_k})}+ \frac{\mu_{k,N_k}}{g_{k,N_k}(u_{k,N_k})}g'_{k,N_k}(u_{k,N_k})\nabla u_{k,N_k}.\label{eqn:gmuik}
\end{equation}
Since
\begin{eqnarray*}
&& g_{k,N_k}(u_{k,N_k})\rightarrow g(u) \text{ uniformly in } D_j^{\circ},\\
&& \frac{g'_{k,N_k}(u_{k,N_k})}{g_{k,N_k}(u_{k,N_k})} \rightarrow \frac{g'(u)}{g(u)} \text{ uniformly in } D_j^{\circ},\\
&& \nabla  \frac{\mu_{k,N_k}}{g_{k,N_k}(u_{k,N_k})} \rightharpoonup \zeta  \text{ weakly in } L^2(D_j^{\circ}),\\
&&\mu_{k,N_k}\rightharpoonup \mu \text{ weakly in } L^2(\Omega_T),\\
&&\nabla u_{k,N_k}\rightarrow \nabla u \text{ strongly in } L^2(\Omega_T),
\end{eqnarray*}
we have, for any $\phi \in L^{\infty}(D_j^{\circ})$,
\begin{eqnarray*}
&&\int_{D_j^{\circ}} \phi \left(g_{k,N_k}(u_{k,N_k})\nabla  \frac{\mu_{k,N_k}}{g_{k,N_k}(u_{k,N_k})}+ \frac{\mu_{k,N_k}}{g_{k,N_k}(u_{k,N_k})}g'_{k,N_k}(u_{k,N_k})\nabla u_{k,N_k}\right) dxdt \\
&&\rightarrow \int_{D_j^{\circ}}\phi \left (g(u)\zeta+\frac{g'(u)}{g(u)}\mu \nabla u\right) dxdt,
\end{eqnarray*}
i.e.
$$
\nabla \mu_{k,N_k}\rightharpoonup \eta\coloneq g(u)\zeta+\frac{g'(u)}{g(u)}\mu \nabla u\ \text{ weakly in } L^1(D_j^{\circ}).
$$
Passing to the limit in \eqref{eqn:gmuk},  we obtain, in  the sense of distribution, that
$$
\eta=-\nabla \Delta u+q''(u)\nabla u+\nabla (-\Delta)^{\frac{1}{2}}u.
$$

Since $q''(u)\nabla u+ \nabla (-\Delta)^{\frac{1}{2}}u \in L^2(\Omega_T)$, we have  $-\nabla \Delta u \in L^1(D_j^{\circ})$, hence
\begin{equation}
\eta=-\nabla \Delta u+q''(u)\nabla u+\nabla (-\Delta)^{\frac{1}{2}}u \text{ a.e. in } D_j^{\circ}
\end{equation}

Since $\frac{1}{g_{k,N_k}(u_{k,N_k})} \rightarrow \frac{1}{g(u)}$ uniformly  in $D_j$, we have
$$
\frac{1}{g_{k,N_k}(u_{k,N_k})}\nabla \mu_{k,N_k} \rightharpoonup \frac{1}{g(u)}\eta \text{ weakly in } L^1(D_j^{\circ}).
$$

Since $\frac{g'_{k,N_k}(u_{k,N_k})}{g^2_{k,N_k}(u_{k,N_k})} \rightarrow \frac{g'(u)}{g^2(u)}$ uniformly  in $D_j$, we have
$$
\frac{g'_{k,N_k}(u_{k,N_k})}{g^2_{k,N_k}(u_{k,N_k})}\mu_{k,N_k}\nabla u_{k,N_k} \rightharpoonup \frac{g'(u)}{g^2(u)}\mu \nabla u \text{ weakly in } L^1(D_j^{\circ}).
$$

Passing to the limit in \eqref{eqn:muikuik}, we have
\begin{eqnarray*}
\zeta&=&\frac{1}{g(u)}\eta-\mu \frac{g'(u)}{g^2(u)}\nabla u\\
&=&\frac{-\nabla \Delta u+q''(u)\nabla u+\nabla (-\Delta)^{\frac{1}{2}}u}{g(u)}-\frac{g'(u)}{g^2(u)}\left(-\Delta u+q'(u)+(-\Delta)^{\frac{1}{2}} u \right)\nabla u
\end{eqnarray*}
on $(D_j)^{\circ}$.
Noticing the value of $\zeta$  on $\Omega_T\backslash D$ doesn't matter since it does not appear on the right hand side of \eqref{eqn:limit}.

{\it{Claim II: For any open set $U\in \Omega_T$ in which  $\nabla \Delta u \in L^p(U)$ for some $p>1$ and $g(u)>0$, we have
\begin{equation}
\zeta=\frac{-\nabla \Delta u+q''(u)\nabla u+\nabla (-\Delta)^{\frac{1}{2}}u}{g(u)}-\frac{g'(u)}{g^2(u)}\left(-\Delta u+q'(u)+(-\Delta)^{\frac{1}{2}} u \right)\nabla u . \label{eqn:zetaexp}
\end{equation} in $U$.}}

To prove this, since
\begin{equation}
\nabla \mu_{k,N_k}=-\nabla \Delta u_{k,N_k}+q''(u_{k,N_k})\nabla u_{k,N_k}+\nabla (-\Delta)^{\frac{1}{2}}u_{k,N_k} \text{ in } \Omega_T \label{eqn:muknk}
\end{equation}
and
\begin{equation}
\nabla \frac{\mu_{k,N_k}}{g_{k,N_k}(u_{k,N_k})} =\frac{1}{g_{k,N_k}(u_{k,N_k})}\nabla \mu_{k,N_k}+\mu_{k,N_k}\cdot \nabla \frac{1}{g_{k,N_k}(u_{k,N_k})} \text{ on } D_j\label{eqn:mugk}.
\end{equation}

The right hand side of \eqref{eqn:muknk} converges weakly to $-\nabla \Delta u+q''(u)\nabla u+\nabla (-\Delta)^{\frac{1}{2}}u$ in $L^q(U)$ for $q=\min\{p,2\}>1$. Hence
$$
\nabla \mu_{k,N_k} \rightharpoonup \eta=-\nabla \Delta u+q''(u)\nabla u+\nabla (-\Delta)^{\frac{1}{2}}u \text{ weakly in } L^q(U).
$$

The right hand side of \eqref{eqn:mugk} converges weakly to
$$\frac{\eta}{g(u)} -\frac{g'(u)}{g^2(u)}\mu\cdot \nabla u$$
in $L^1(U\cap D_j)$ for each $j$ and therefore $$\zeta=\frac{-\nabla \Delta u+q''(u)\nabla u+\nabla (-\Delta)^{\frac{1}{2}}u}{g(u)}-\frac{g'(u)}{g^2(u)}\left(-\Delta u+q'(u)+(-\Delta)^{\frac{1}{2}} u \right)\nabla u $$ a.e. in $ U\cap D$.  and the definition of $\zeta$ can be extended to $U\backslash D$ by our integrability assumption on $u$. Define
$$
\tilde \Omega_T=\{U\subset \Omega_T:\nabla \Delta u \in L^p(U) \text{ for some } p>1; g(u)>0 \text{ on } U \text{ depending on } U\}.
$$
Then $\tilde \Omega_T$ is open and $\zeta$ is defined by  \eqref{eqn:zetaexp} on  $\tilde \Omega_T$. Since $|\Omega_T\backslash B|=0$ ,  $M(u)=0$ on $\Omega_T\backslash P$ and
$$
\Omega_T\backslash \{D\cup \tilde \Omega_T\}\subset \{\Omega_T\backslash B\} \cup \{\Omega_T\backslash P\},
$$
we can take the value  of $\zeta$ to be zero outside $D\cup \Omega_T$, sand it won't affect the integral on the right side of \eqref{eqn:dg}.

Lastly the energy inequality \eqref{eqn:ineq}  follows by taking limit in  the energy inequality for $u_{k,N_k}$.

This finishes the proof of Theorem \ref{thm:dg}.

\section{Simulations}\label{sec:simulation}

In this section, we use the proposed phase field model  to simulate the climb motions of prismatic dislocation loops, incorporating the conservative motion and nonconservative motion. We use the evolution equation in Eqs.~\eqref{eqn:ch} without the factor $g(u)$ on the right-hand side, i.e.,
\begin{flalign}\label{eqn:model000}
 \partial_t u+\beta \mu=\nabla \cdot \left(M(u)\nabla \frac{\mu}{g(u)}\right),
\end{flalign}
together with Eqs.~\eqref{eqn:chem} and \eqref{equ:clf}.
 Recall that the nonconservative climb motion will result into the shrinking and growing of the dislocation loops \cite{Hirth-Lothe}, whereas the  self-climb is a conservative motion, which will keep the enclosed area of a prismatic loop unchanged \cite{Kroupa1961,Niu2017,Niu2019}.

 In the simulations, we choose the simulation domain $ \Omega = [-\pi, \pi]^2$ and mesh size $dx=dy=2\pi/M$ with $M=64$. Periodic boundary conditions are used for the simulation domain. The small parameter in the phase field model $\varepsilon=dx$. The simulation domain corresponds to a physical domain of size $(300b)^2$, i.e., $b=2\pi/300$. Under this setting, the parameter $H_0$ in the phase field model calculated in the paper \cite{Niu2021} 
is $H_0=52.65\left(2( 1-\nu ) /\mu b^2\right)$. The prismatic loops are in the counterclockwise direction meaning vacancy loops, unless otherwise specified.

In the numerical simulations, we use the pseudospectral method: All the spatial partial derivatives are calculated in the Fourier space using FFT. For the time discretization, we use  the forward Euler method. The climb force generated by dislocations $f_{\rm{cl}}^{\rm d}$ is calculated by FFT using Eq.~\eqref{equ:clf}. We regularize the function $g(u)$ in the denominator in Eq.~\eqref{eqn:model000} as $\sqrt{g(u)^2+e_0^2}$ with small parameter $e_0=0.005$. In the initial configuration of a simulation, $\phi$ in the dislocation core region is set to be a $\tanh$ function with width $3\varepsilon$. The location of the dislocation loop is identified by the contour line of $u=0$.


\subsection{Evolution of an elliptic prismatic loop under the combined climb effect}

\begin{figure}[h]
\centering
\subfigure[\label{fig.elliptic.with.selfclimb}]{\includegraphics[width=3.in]{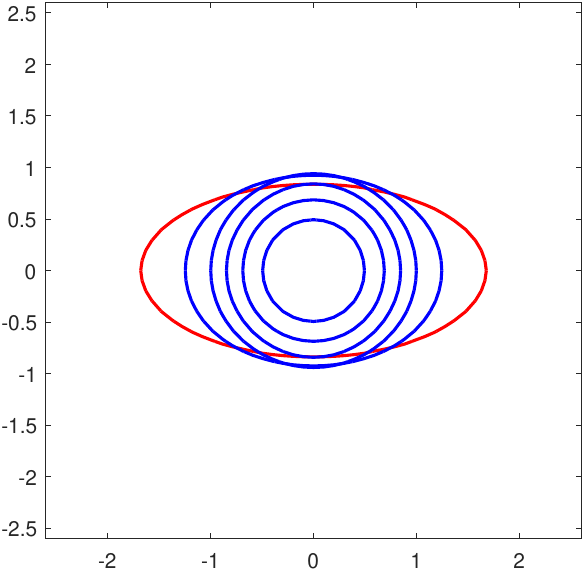}}
\subfigure[\label{fig.elliptic.without.selfclimb}]{\includegraphics[width=3.in]{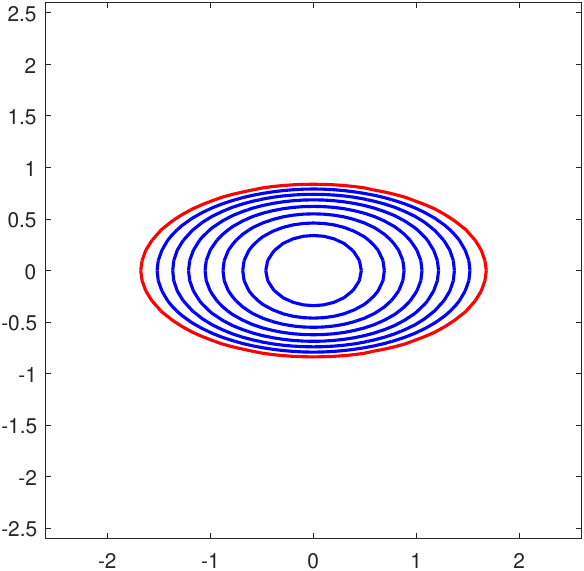}}
\label{fig.elliptic}
\caption{shrinking  of an elliptic prismatic loops by climb with/without self-climb.}
\end{figure}

\begin{figure}[h]
\centering
\subfigure{\includegraphics[width=3.in]{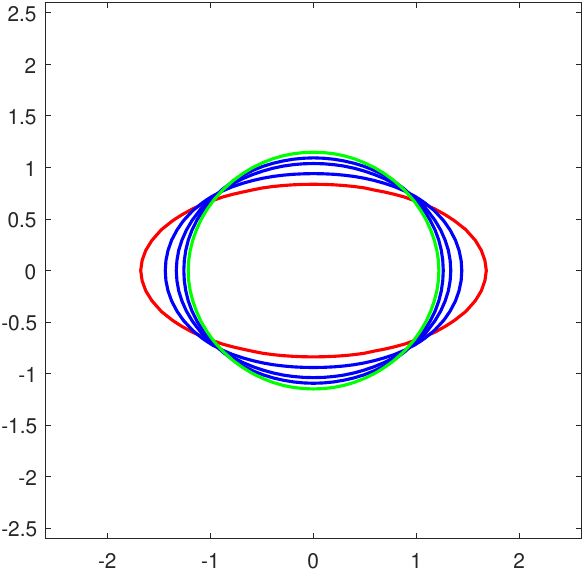}}
\caption{Evolution of an elliptic prismatic loop  only by self-climb effect using the phase field model. Red ellipse is the initial state, and green circle is the final state.}
\label{fig.elliptic.only.selfclimb}
\end{figure}

In the first numerical example, we simulate evolution of an elliptic prismatic loop using the phase field model, see Fig.~\ref{fig.elliptic} and Fig.~\ref{fig.elliptic.only.selfclimb}. The two axes of the initial elliptic profile are $l_1=80b$ and $l_2=40b$. Fig.~\ref{fig.elliptic.with.selfclimb} shows the elliptic prismatic loop will not directly shrink, due to the self-climb effect, and there is a trend to evolve to a circle in the shrinking process. Fig.~\ref{fig.elliptic.without.selfclimb} shows  that without the self-climb effect, the elliptic loop directly shrink until vanishing. The time of shrinking of loop with self-climb is much more than the time without self-climb. The shapes  are also {  totally different}  in the process. These will influence the pattern of the interactions of two loops, see details in the simulations; see Sec.(\ref{simulation.twoloops}). Moreover, we show the evolution of an elliptic prismatic loop only by self-climb  using the phase field model, seeing Fig.~\ref{fig.elliptic.only.selfclimb}, to illustrate the effect of the self-climb effect. Red ellipse is the initial state,  and the loop converges to the equilibrium shape of a circle (green circle) under its self-stress. The area enclosed by a prismatic loop is conserved during the self-climb motion. More simulation information  about the self-climb effect can be found in our previous papers \cite{Niu2017, Niu2019, Niu2021}. 

 \subsection{The interactions between the circular  prismatic loops under the combined climb effect}\label{simulation.twoloops}

\begin{figure}[htbp]
\centering
\subfigure[\label{fig.twoloops-with-selfclimb-1}]{\includegraphics[width=2.in]{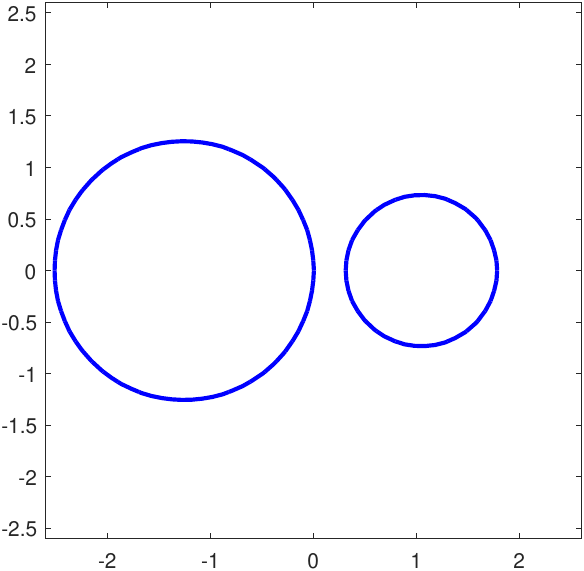}}
\subfigure[\label{fig.twoloops-with-selfclimb-2}]{\includegraphics[width=2.in]{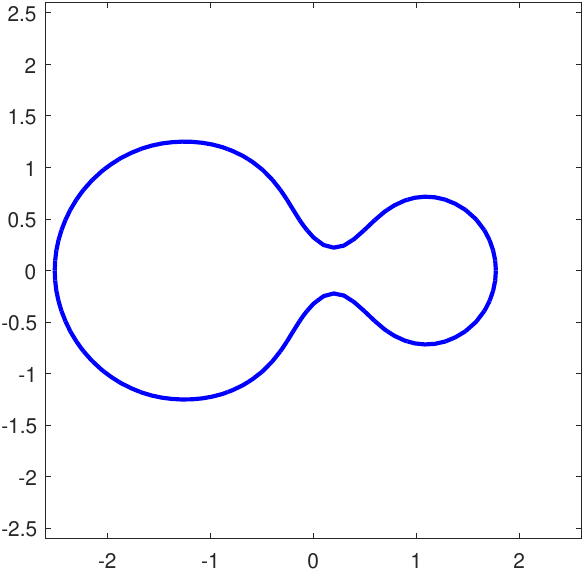}}
\subfigure[\label{fig.twoloops-with-selfclimb-3}]{\includegraphics[width=2.in]{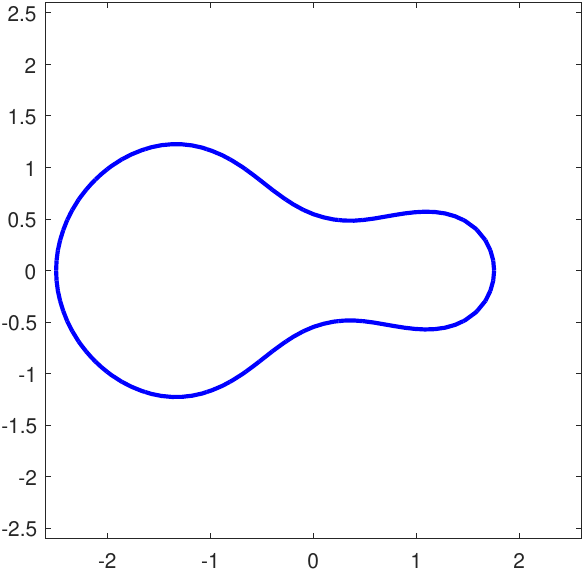}}
\subfigure[\label{fig.twoloops-with-selfclimb-4}]{\includegraphics[width=2.in]{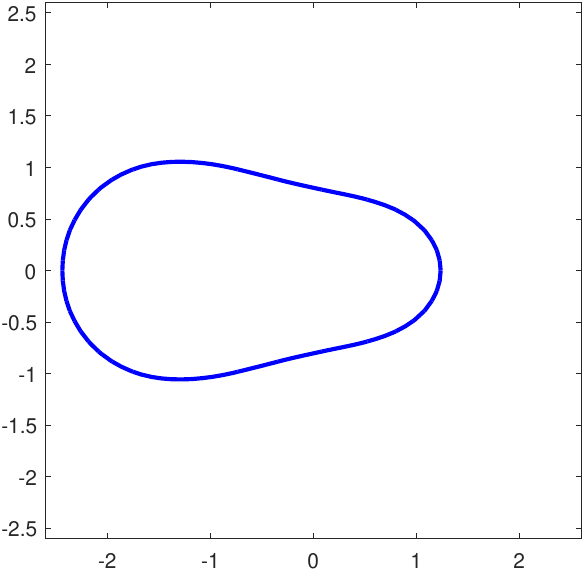}}
\subfigure[\label{fig.twoloops-with-selfclimb-5}]{\includegraphics[width=2.in]{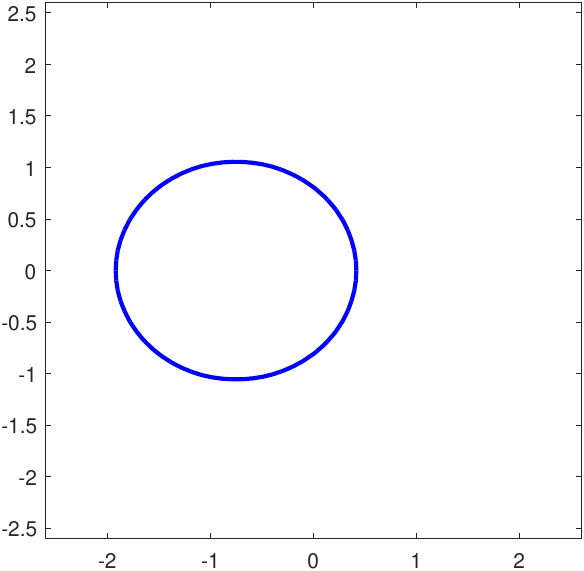}}
\subfigure[\label{fig.twoloops-with-selfclimb-6}]{\includegraphics[width=2.in]{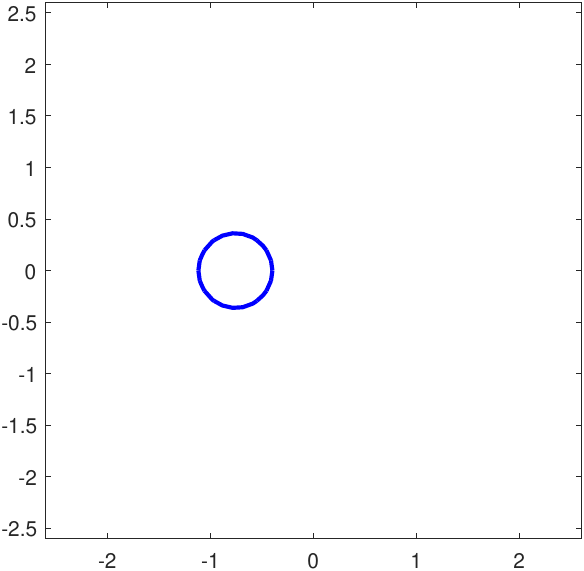}}
	\caption{The interaction of two circular prismatic loops by climb with self-climb.}\label{fig.twoloops-with-selfclimb}
\end{figure}

\begin{figure}[htbp]
\centering
\subfigure[\label{fig.twoloops-without-selfclimb-1}]{\includegraphics[width=2.in]{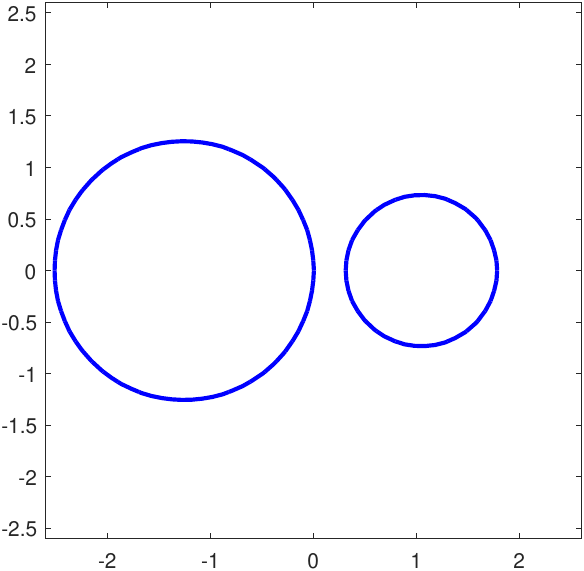}}
\subfigure[\label{fig.twoloops-without-selfclimb-2}]{\includegraphics[width=2.in]{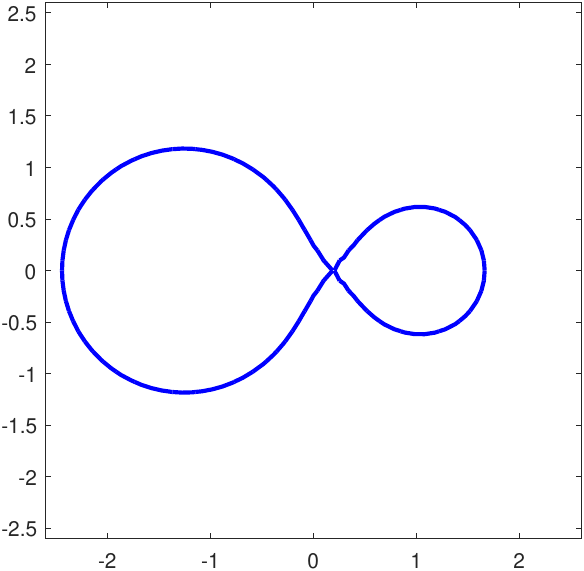}}
\subfigure[\label{fig.twoloops-without-selfclimb-3}]{\includegraphics[width=2.in]{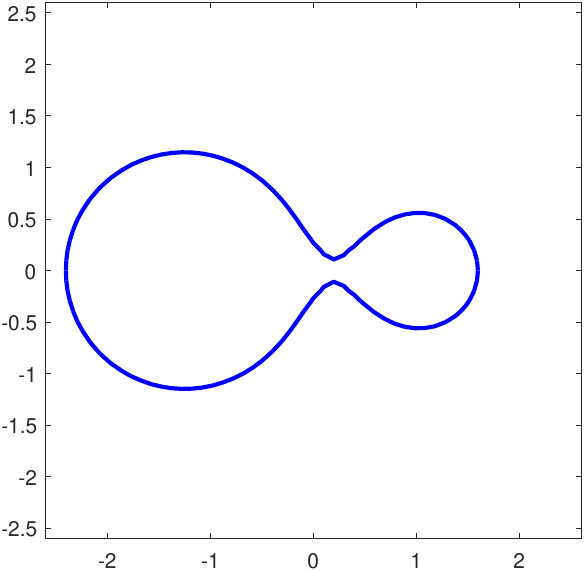}}
\subfigure[\label{fig.twoloops-without-selfclimb-4}]{\includegraphics[width=2.in]{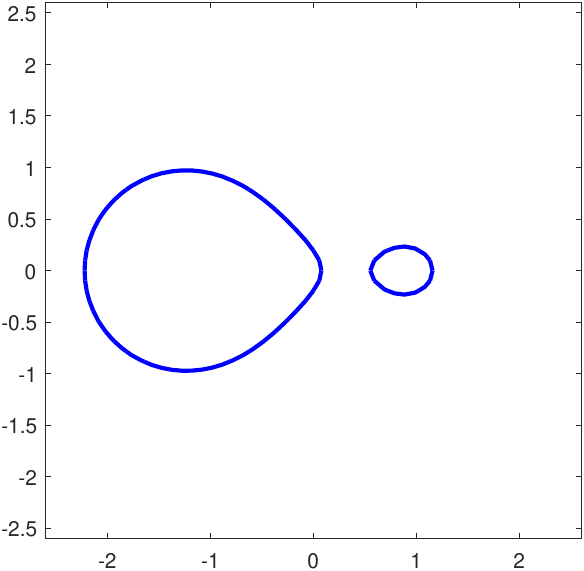}}
\subfigure[\label{fig.twoloops-without-selfclimb-5}]{\includegraphics[width=2.in]{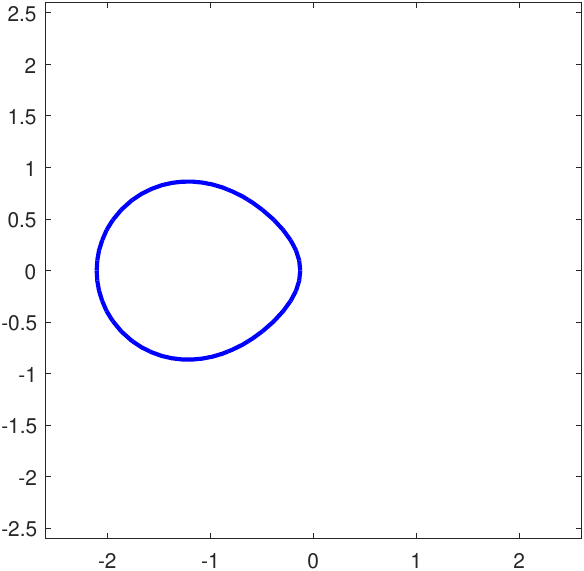}}
\subfigure[\label{fig.twoloops-without-selfclimb-6}]{\includegraphics[width=2.in]{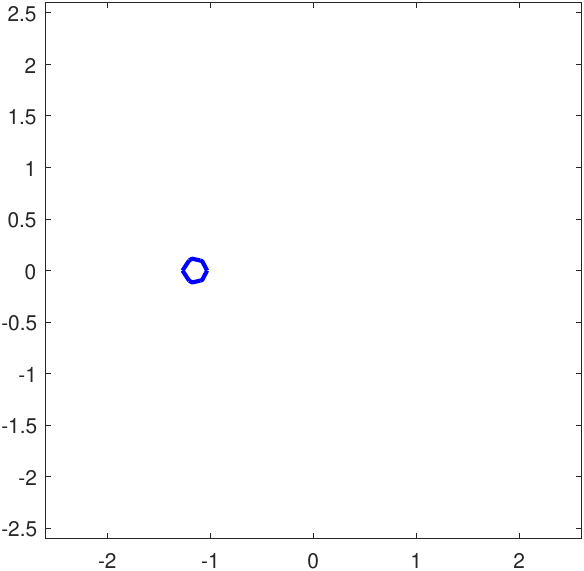}}
	\caption{The interaction of two circular prismatic loops by climb without self-climb.}\label{fig.twoloops-without-selfclimb}
\end{figure}

\begin{figure}[htbp]
\centering
\subfigure[]{\includegraphics[width=2.in]{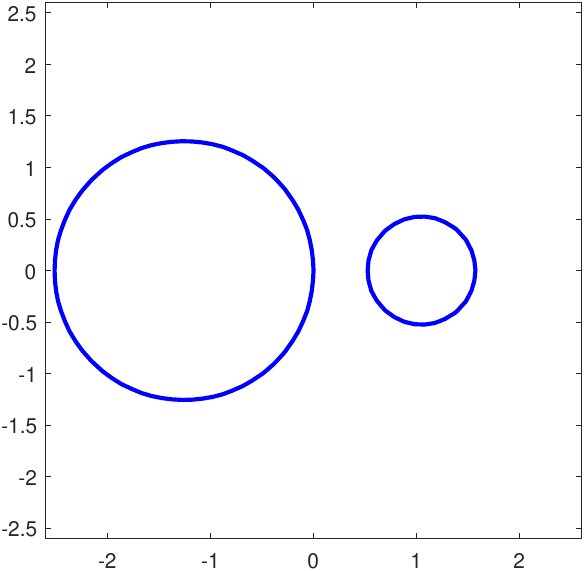}}
\subfigure[]{\includegraphics[width=2.in]{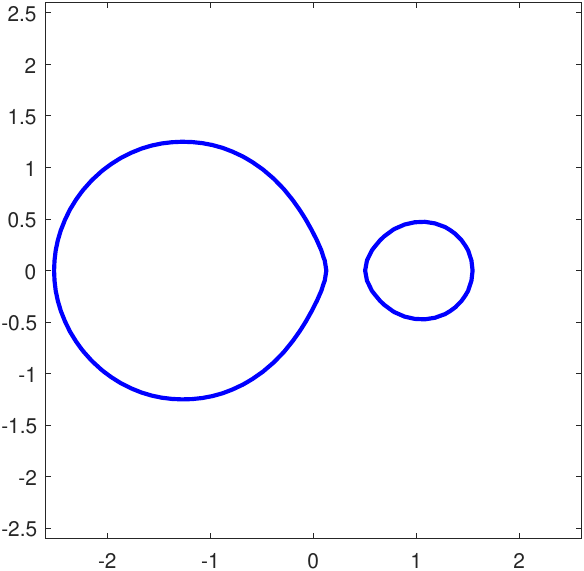}}
\subfigure[]{\includegraphics[width=2.in]{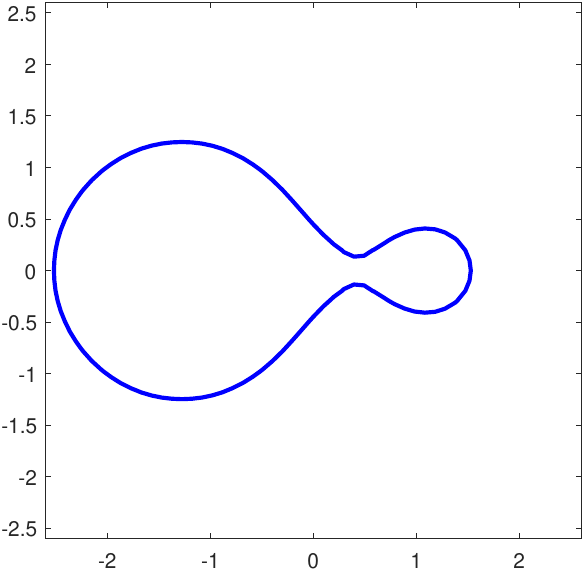}}
\subfigure[]{\includegraphics[width=2.in]{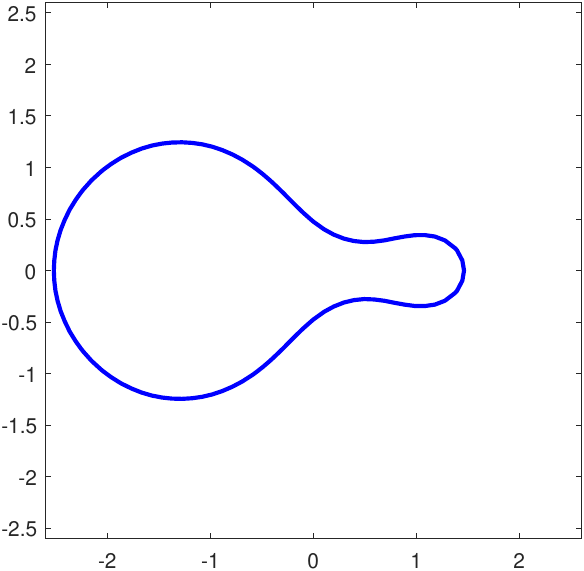}}
\subfigure[]{\includegraphics[width=2.in]{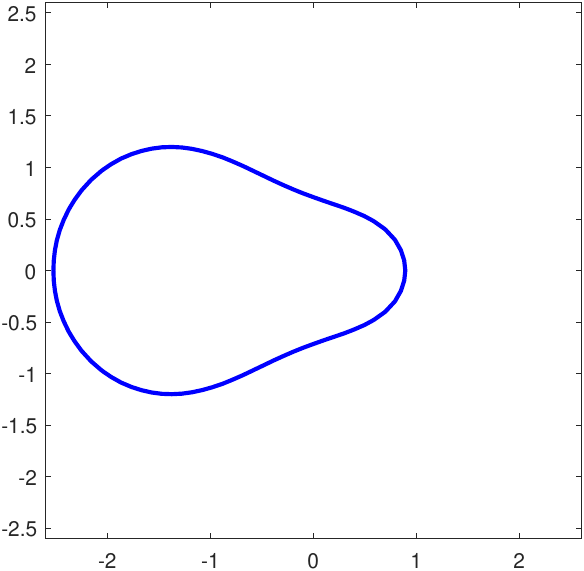}}
\subfigure[]{\includegraphics[width=2.in]{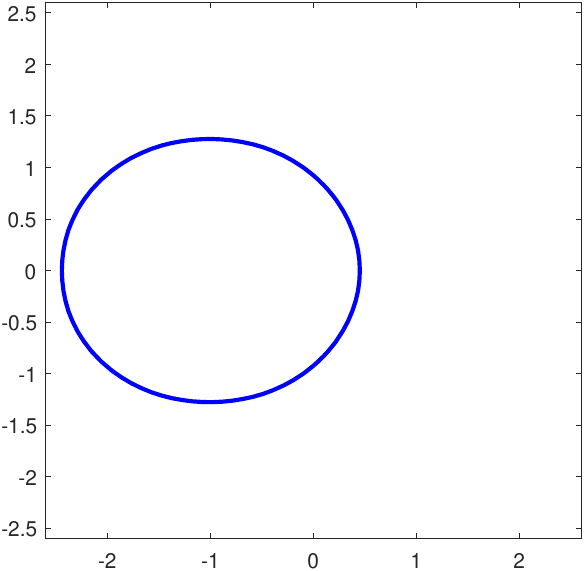}}
	\caption{Coalescence of two prismatic loops only by self-climb under their elastic interaction obtained by the phase field model.}\label{fig.twoloops-only-selfclimb}
\end{figure}

  In this subsection, we use our  phase field model to simulate  the interaction of two circular prismatic loops for three conditions: with self-climb, without self-climb and only with self-climb. The detailed shrinking process obtained by our simulations are shown in Fig.~\ref{fig.twoloops-with-selfclimb}, Fig.~\ref{fig.twoloops-without-selfclimb} and  Fig.~\ref{fig.twoloops-only-selfclimb}. The two loops are attracted to each other by self-climb under the elastic interaction between them for all these three conditions, but the later change of the shapes are totally different. For the simulation of dislocation climb with the self-climb effect, firstly, the two loops are attracted to each other by self-climb. When the two loops meet, they quickly combine into a single loop; see Fig.~\ref{fig.twoloops-with-selfclimb}(a-b). The combined single loop eventually evolves into a circular shape; see Fig.~\ref{fig.twoloops-with-selfclimb}(c)-(e). Finally the circular loops shrink and vanish; see Fig.~\ref{fig.twoloops-with-selfclimb}(f).  For the simulation of dislocation climb without the self-climb effect, see Fig.~\ref{fig.twoloops-without-selfclimb}. Firstly, the two loops are attracted to each other under self-stress; see Fig.~\ref{fig.twoloops-without-selfclimb}(a)-(c), but quickly they separate due to the non-conservative climb effect; see Fig.~\ref{fig.twoloops-without-selfclimb}(d). The small loop  {  vanishes first} in the shrinking process;  see Fig.~\ref{fig.twoloops-without-selfclimb}(e). Finally the larger  loop shrinks and vanishes; see Fig.~\ref{fig.twoloops-without-selfclimb}(f). Comparing these two climb interaction processes  with and without self-climb effect, we {  conclude that}  even though {  both}  loops will vanish {  eventually}, the processes {  are quite different}. With the effect of self-climb, these two close loop will coalescence  first when they shrink. Without the self-climb, these two loops  will shrink directly and simply after the quick connecting and separation. The total time for the shrinking of these two loops {  differs} greatly.  It takes longer time for the loops {  to shrink} with the self-climb effect than without the self-climb effect. Fig.~\ref{fig.twoloops-with-selfclimb} and Fig.~\ref{fig.twoloops-without-selfclimb} give details of the patterns in these two shrinking process and show the great difference, which will help us to understand the formation process of the patterns and predict the stable state of the patterns in the physics experiments. Moreover, for understanding the self-climb effect in the interactions of the two loop, we show the detailed coalescence process only by self-climb in Fig.~\ref{fig.twoloops-only-selfclimb}. Firstly, the two loops are attracted to each other by self-climb under the elastic interaction. They quickly combine into a single loop after meeting; see Fig.~\ref{fig.twoloops-only-selfclimb}(a-c). The combined single loop eventually evolves into a stable, circular shape; see Fig.~\ref{fig.twoloops-only-selfclimb}(d)-(f). It is noteworthy that the area of the final circle are equal to the total area of the initial two circles  theoretically, and  these two areas also agree well in numerical simulation. More simulation information  about the self-climb effect of the interactions of loops can be found in our previous papers \cite{Niu2017, Niu2019, Niu2021}.

\section{Conclusions and discussion}\label{sec:conclusions}


 We have performed asymptotic analysis to show  that our phase field model yields the correct climb velocity  in the sharp interface limit including the self-climb contribution.

 we have validated our phase field model by numerical simulations and compared the evolutions of the elliptic loops with and without self-climb, the interactions of two loops with and without self-climb. The simulation results show the  {  big }  difference {  in } the evolution time and the pattern with and without self-climb contribution.

%

Self-climb by vacancy pipe diffusion is the dominant dislocation climb mechanism at not very high temperature in irradiated materials. At high temperature, dislocation climb by vacancy bulk diffusion also becomes important, the contribution portion of both climb motions can be adjusted by the parameter
$\beta$ depending on the physical material and situation, and these increase the applicability of the phase field model in physics. This phase field model combines these two climb motions in a single evolution equation, which can simulate the combined climb motion of interactions of many loops. This model {  can be  easily}  generalized to many loops and the interactions of loops in 3-dimension space. It also provides a convenient base to simulate dislocation climb-glide motion.

\textbf{ACKNOWLEDGEMENTS}

 Both authors thank Prof. Yang Xiang for the helpful discussions. X.H. Niu's research is supported by National Natural Science Foundation of China under the grant number 11801214 and the Natural Science Foundation of Fujian Province of China under the grant number 2021J011193.
 X. Yan's research is supported by a Research Excellence Grant, CLAS Dean's Summer Research Grant from University of Connecticut and Simons travel grant \#947054.

\FloatBarrier
\bibliographystyle{siam}
\bibliography{references}

\end{document}